\documentclass[preprint,12pt]{elsarticle}

\usepackage{amsmath}    
\usepackage{graphicx,amssymb}
\usepackage{epsfig}

\usepackage{subcaption}     

\allowdisplaybreaks[1]

\begin{document}

\def\eR{\mathbb{R}}     

\title{Fractal properties of Bessel functions}

\author[fer]{L.\ Korkut}
\ead{luka.korkut@fer.hr}

\author[fer]{D.\ Vlah\corref{corresp}}
\ead{domagoj.vlah@fer.hr}

\author[fer]{V.\ \v Zupanovi\'c}
\ead{vesna.zupanovic@fer.hr}


\address[fer]{University of Zagreb, Faculty of Electrical Engineering and Computing, Unska 3, 10000 Zagreb, Croatia}

\journal{arXiv.org}

\begin{abstract}
A fractal oscillatority of solutions of second-order differential equations near infinity is measured by oscillatory and phase dimensions. The phase dimension is defined as a box dimension of the trajectory $(x,\dot{x})$ in $\eR^2$ of a solution $x=x(t)$, assuming that $(x,\dot{x})$ is a spiral converging to the origin. In this work, we study the phase dimension of the class of second-order nonautonomous differential equations with oscillatory solutions including the Bessel equation. We prove that the phase dimension of Bessel functions is equal to $4/3$, and that the corresponding trajectory is a wavy spiral, exhibiting an interesting behavior. The phase dimension of a generalization of the Bessel equation has been also computed.
\end{abstract}
\medskip

\date{}

\maketitle



\newtheorem{theorem}{Theorem}
\newtheorem{cor}{Corollary}
\newtheorem{prop}{Proposition}
\newtheorem{lemma}{Lemma}
\newdefinition{defin}{Definition}
\newdefinition{remark}{Remark}

\font\csc=cmcsc10

\def\esssup{\mathop{\rm ess\,sup}}
\def\essinf{\mathop{\rm ess\,inf}}
\def\wo#1#2#3{W^{#1,#2}_0(#3)}
\def\w#1#2#3{W^{#1,#2}(#3)}
\def\wloc#1#2#3{W_{\scriptstyle loc}^{#1,#2}(#3)}
\def\osc{\mathop{\rm osc}}
\def\var{\mathop{\rm Var}}
\def\supp{\mathop{\rm supp}}
\def\Cap{{\rm Cap}}
\def\norma#1#2{\|#1\|_{#2}}

\def\C{\Gamma}

\let\text=\mbox

\catcode`\@=11
\let\ced=\c
\def\a{\alpha}
\def\b{\beta}
\def\d{\delta}
\def\g{\lambda}
\def\o{\omega}
\def\q{\quad}
\def\n{\nabla}
\def\s{\sigma}
\def\div{\mathop{\rm div}}
\def\sing{{\rm Sing}\,}
\def\singg{{\rm Sing}_\ty\,}

\def\A{{\cal A}}
\def\F{{\cal F}}
\def\H{{\cal H}}
\def\W{{\bf W}}
\def\M{{\cal M}}
\def\N{{\cal N}}
\def\S{{\cal S}}

\def\ty{\infty}
\def\e{\varepsilon}
\def\f{\varphi}
\def\:{{\penalty10000\hbox{\kern1mm\rm:\kern1mm}\penalty10000}}
\def\ov#1{\overline{#1}}
\def\D{\Delta}
\def\O{\Omega}
\def\pa{\partial}

\def\st{\subset}
\def\stq{\subseteq}
\def\pd#1#2{\frac{\pa#1}{\pa#2}}
\def\sgn{{\rm sgn}\,}
\def\sp#1#2{\langle#1,#2\rangle}

\newcount\br@j
\br@j=0
\def\q{\quad}
\def\gg #1#2{\hat G_{#1}#2(x)}
\def\inty{\int_0^{\ty}}
\def\od#1#2{\frac{d#1}{d#2}}

\def\bg{\begin}
\def\eq{equation}
\def\bgeq{\bg{\eq}}
\def\endeq{\end{\eq}}
\def\bgeqnn{\bg{eqnarray*}}
\def\endeqnn{\end{eqnarray*}}
\def\bgeqn{\bg{eqnarray}}
\def\endeqn{\end{eqnarray}}

\def\bgeqq#1#2{\bgeqn\label{#1} #2\left\{\begin{array}{ll}}
\def\endeqq{\end{array}\right.\endeqn}

\def\abstract{\bgroup\leftskip=2\parindent\rightskip=2\parindent
        \noindent{\bf Abstract.\enspace}}
\def\endabstract{\par\egroup}

\def\udesno#1{\unskip\nobreak\hfil\penalty50\hskip1em\hbox{}
             \nobreak\hfil{#1\unskip\ignorespaces}
                 \parfillskip=\z@ \finalhyphendemerits=\z@\par
                 \parfillskip=0pt plus 1fil}
\catcode`\@=11

\def\cal{\mathcal}

\def\eN{\mathbb{N}}
\def\Ze{\mathbb{Z}}
\def\Qu{\mathbb{Q}}
\def\Ce{\mathbb{C}}

\def\osd{\mathrm{osd}\,}

\section{Introduction and motivation}

The fractal oscillatority of solutions of different types od second-order linear differential equations has been recently considered by Kwong, Pa{\v{s}}i{\'c}, Tanaka and Wong. The Euler type equation has been studied in \cite{pasiceuler} and  \cite{pasic}, the Hartman-Wintner type equations in \cite{kpw}, half linear equations in \cite{pasicwo}, and finally the Bessel equation in \cite{mersat}. In all of this work the fractal oscillatority is considered in the sense of the oscillatory dimension, at first introduced in Pa\v si\'c, \v Zubrini\'c and \v Zupanovi\'c \cite{chirp}. The oscillatory dimension of a solution $x(t)$ is defined as the box dimension of a graph of the function $X(\tau):=x(\frac{1}{\tau})$ near $\tau=0$.

On the other hand, the fractal properties of spiral trajectories of dynamical systems in the phase plane have been studied by \v Zubrini\'c and \v Zupanovi\'c, see e.g.\ \cite{zuzu} and \cite{zuzulien}. From their work the concept of the phase dimension has arisen and has finally been introduced in \cite{chirp}. They adapted standard idea of phase plane analysis to fractal analysis of solutions of second-order nonlinear autonomous differential equations.

These ideas motivated us to study a fractal connection between the oscillatory and phase dimensions for a class of oscillatory functions, see \cite{cswavy}. In that study, we discovered a specific type of spirals with a nondecreasing radius function, related to chirp-like solutions of a class of equations considered by Kwong, Pa{\v{s}}i{\'c}, Tanaka and Wong, which we call the wavy spirals.

A model for chirp-like behavior of solutions developed in \cite{cswavy} could not handle the specific behavior of Bessel functions, so to determine the phase dimension of Bessel functions we generalize the technique used in \cite{cswavy}. Notice that the oscillatory dimension of Bessel functions was already considered in \cite{mersat}.

In this work we actually consider some generalization of the Bessel equation that is motivated by the generalization introduced in \cite{mersat}, for which we determine the phase dimension in our main result, Theorem \ref{tm_Bessel}. The standard Bessel equation now becomes a specific case of this generalization. In order to prove Theorem \ref{tm_Bessel}, we first obtain a new version of some theorems from \cite{zuzu}.

We discover an interesting property related to wavy spirals. By varying values of parameters in our generalized Bessel equation, we go from spirals with no waves to spirals with ``big" waves.
\smallskip

\section{Definitions}

We first introduce some definitions and notation.
For $A\st\eR^N$  bounded
we define the \emph{$\e$-neighborhood} of $A$ by $A_\e:=\{y\in\eR^N\:d(y,A)<\e\}$.
By the {\it lower $s$-dimensional  Minkowski content} of $A$, $s\ge0$ we mean
$$
\M_*^s(A):=\liminf_{\e\to0}\frac{|A_\e|}{\e^{N-s}},
$$
and analogously for the {\it upper $s$-dimensional Minkowski content} $\M^{*s}(A)$.
Now we can introduce the lower and upper box dimensions of $A$ by
$$
\underline\dim_BA:=\inf\{s\ge0\:\M_*^s(A)=0\}
$$
 and analogously
$\ov\dim_BA:=\inf\{s\ge0\:\M^{*s}(A)=0\}$.
If these two values coincide, we call it simply the box dimension of $A$, and denote by $\dim_BA$.

For more details on these definitions see e.g.\ Falconer \cite{falc}, \cite{zuzu} and \cite{cras}.

\smallskip

Assume now that $x$ is of class $C^1$ and $t_0>0$. We say that $x$ is a {\it phase oscillatory} function if the following condition holds: the set
$\C=\{(x(t),\dot x(t)):t\in[t_0,\ty)\}$ in the plane is a spiral converging to the origin.

By the \emph{spiral} here we mean the graph of a function $r=f(\f)$, $\f\geq\f_1>0$, in polar coordinates, where
\begin{equation}\label{def_spiral}
\left\{\begin{array}{l}
f:[\f_1,\infty)\rightarrow(0,\infty) \textrm{ is such that } f(\f)\to 0 \textrm{ as } \f\to\infty,\\
f \textrm{ is \emph{radially decreasing} (i.e., for any fixed } \f\geq\f_1\\
\textrm{the function } \mathbb{N}\ni k\mapsto f(\f+2k\pi) \textrm{ is decreasing)} .
\end{array}\right.
\end{equation}
Depending on the context, by the spiral here we also mean the graph of a function $r=g(\f)$, $\f\leq\f'_1<0$, in polar coordinates, where for $h(\f)=g(-\f)$, $\forall \f\geq|\f'_1|$, the graph of a function $r=h(\f)$, $\f\geq|\f'_1|>0$, given in polar coordinates, satisfies (\ref{def_spiral}). It is easy to see that the spiral given by function $g$ is a mirror image of the spiral given by function $h$, regarding $x$-axis. We also say that a graph of function $r=f(\f)$, $\f\geq\f_1>0$, in polar coordinates, is a \emph{spiral near the origin} if there exists $\f_2\geq\f_1$ such that a graph of function $r=f(\f)$, $\f\geq\f_2$ is a spiral.

The {\it phase dimension} $\dim_{ph}(x)$ of function $x:[t_0,\infty)\to\eR$ of class $C^1$ is defined as the box dimension of the corresponding planar curve $\C=\{(x(t),\dot x(t)):t\in[t_0,\ty)\}$.

\smallskip

We use a result for the box dimension of spiral $\C$ defined by $r=\f^{-\a}$, $\f\ge\f_0>0$,  $\dim_B\C=2/(1+\a)$ when $0<\a\le1$, see Tricot \cite[p.\ 121]{tricot} and some generalizations from \cite{zuzu}.

The phase dimension is a fractal dimension, introduced in the study of chirp-like solutions of second order ODEs, see \cite{chirp}. Fractal dimensions are a well known tool in study of dynamics, see \cite{fdd}.

\smallskip

For two real functions $f(t)$ and $g(t)$ of a real variable we write $f(t)\simeq g(t)$ as $t\to0$ (as $t\to\ty$) if there exist positive constants $C$ and $D$ such that $C\,f(t)\le g(t)\le D\,f(t)$, for all $t$ sufficiently large. For example, for a function $F:U\to V$ with $U,V\st\eR^2$, $V=F(U)$, the condition $|F(x_1)-F(x_2)|\simeq
|x_1-x_2|$ means that $f$ is a bi-Lipschitz mapping, i.e.,\ both $F$ and $F^{-1}$ are Lipschitzian.

\smallskip

We write $f(t)\sim g(t)$ if $f(t)/g(t)\to1$ as $t\to\ty$. Also, if $k$  is a fixed positive integer, for two functions $f$ and $g$ of class $C^k$ we write,
$$
f(t)\sim_k g(t) \ \mbox{as}\ t\to\ty ,
$$
if $f^{(j)}(t)\sim g^{(j)}(t)$ as $t\to\ty$ for all $j=0,1,...,k$.

For example, $\frac{(t-1)^{4-\a}}{t^4}\sim_3t^{-\a}$ as $t\to\ty$, for $\a\in(0,1)$.

\smallskip

We write $f(t)=O(g(t))$ as $t\to\ty$ if there exists a positive constant $C$ such that $|f(t)|\leq C|g(t)|$ for all $t$ sufficiently large. Similary, we write $f(t)=o(g(t))$ as $t\to\ty$ if for every positive constant $\varepsilon$ it holds $|f(t)|\leq \varepsilon|g(t)|$ for all $t$ sufficiently large.

\section{Phase dimension of Bessel functions}

From the point of view of a fractal geometry we study a spiral generated by Bessel functions and generalized Bessel functions.

The \emph{Bessel equation of order $\nu$}, widely known in literature, see e.g. \cite[p.~98]{lebedev}, is the linear second-order ordinary differential equation given by
\begin{equation} \label{Bessel_equation}
t^2 x''(t)+t x'(t)+(t^2-\nu^2)x(t)=0 ,
\end{equation}
with one parameter $\nu\in\eR$. It has two linearly independent solutions, known as \emph{Bessel functions} of the first and second kind of order $\nu$, designated $J_{\nu}$ and $Y_{\nu}$, respectively.

Pa{\v{s}}i{\'c} and Tanaka in \cite{mersat} studied the oscillatory dimension of the Bessel equation and one specific generalization of this equation that they introduced. Motivated by their generalization we introduce this class of equations
\begin{equation} \label{gen_Bessel_equation}
t^2 x''(t)+t(2-\mu)x'(t)+(t^2-\nu^2)x(t)=0 ,
\end{equation}
where $\nu,\mu\in\eR$ are parameters. By setting $\mu=1$ we get the standard Bessel equation (\ref{Bessel_equation}). Two linearly independent solutions to equation (\ref{gen_Bessel_equation}) define functions $\widetilde{J}_{\nu,\mu}$ and $\widetilde{Y}_{\nu,\mu}$, which we call generalized Bessel functions. By solving (\ref{gen_Bessel_equation}), we see that $\widetilde{J}_{\nu,\mu}$ and $\widetilde{Y}_{\nu,\mu}$ can be written in terms of $J_{\widetilde{\nu}}$ and $Y_{\widetilde{\nu}}$, respectively, where $\widetilde{\nu}=\sqrt{\left(\frac{\mu-1}{2}\right)^2+\nu^2}$. We get
\begin{eqnarray} \label{gen_Bessel_functions}
\widetilde{J}_{\nu,\mu}(t) & = & t^{\frac{\mu-1}{2}} J_{\widetilde{\nu}}(t) , \nonumber\\
\widetilde{Y}_{\nu,\mu}(t) & = & t^{\frac{\mu-1}{2}} Y_{\widetilde{\nu}}(t) . \nonumber
\end{eqnarray}

As $J_{\nu}(t)=O\left(t^{-\frac{1}{2}}\right)$ and $Y_{\nu}(t)=O\left(t^{-\frac{1}{2}}\right)$, for every $\nu\in\eR$, it follows that $\widetilde{J}_{\nu,\mu}(t) = O\left(t^{\frac{\mu-2}{2}}\right)$ and $\widetilde{Y}_{\nu,\mu}(t) = O\left(t^{\frac{\mu-2}{2}}\right)$, for every $\nu,\mu\in\eR$. Now, assume that $\tau_0>0$. It is easy to see the following:
\begin{itemize}
\item
If $\mu>2$ then $\widetilde{J}_{\nu,\mu}(t)$ and $\widetilde{Y}_{\nu,\mu}(t)$ are unbounded for $t\in[\tau_0,\infty)$. Trivially, $\widetilde{J}_{\nu,\mu}(t)$ and $\widetilde{Y}_{\nu,\mu}(t)$ are not phase oscillatory functions.
\item
If $\mu=2$ then $\widetilde{J}_{\nu,\mu}(t)$ and $\widetilde{Y}_{\nu,\mu}(t)$ are bounded for $t\in[\tau_0,\infty)$, but $\lim_{t\to\infty} \widetilde{J}_{\nu,\mu}(t)$ and $\lim_{t\to\infty} \widetilde{Y}_{\nu,\mu}(t)$ do not exist, so $\widetilde{J}_{\nu,\mu}(t)$ and $\widetilde{Y}_{\nu,\mu}(t)$ are not phase oscillatory functions.
\item
If $\mu<0$ then $\widetilde{J}_{\nu,\mu}(t)$ and $\widetilde{Y}_{\nu,\mu}(t)$ are phase oscillatory functions, but with a trivial phase dimension, $\dim_{ph}\left(\widetilde{J}_{\nu,\mu}\right)=\dim_{ph}\left(\widetilde{J}_{\nu,\mu}\right)=1$, for every $\nu\in\eR$.
\end{itemize}
The remaining case $\mu\in(0,2)$, which will produce phase oscillatory functions with a nontrivial phase dimension, we describe in Theorem \ref{tm_Bessel}.

\begin{theorem}{\label{tm_Bessel}} {\rm (Generalized Bessel spiral)}
Let $\mu\in(0,2)$, $\nu\in\eR$ and $\tau_0>0$. Let $x(t)=\widetilde{J}_{\nu,\mu}(t)$, for $t\in[\tau_0,\infty)$, and let continuous function $\f(t)$ be given by $\tan \f(t)=\frac{\dot x(t)}{x(t)}$, for $t\in[\tau_0,\infty)$.

Then the planar curve $\Gamma=\{(x(t),\dot x(t)): t\in[\tau_0,\ty)\}$ is a spiral near the origin $r=g(\f)$, $\f\in(-\ty,-\phi_0]$, satisfying $g(\f)\simeq |\f|^{-\frac{2-\mu}{2}}$, in polar coordinates. Phase dimension
\begin{equation} \label{dim_rezultat_tm}
\dim_{ph}(x)=\dim_B\Gamma=\frac{4}{4-\mu} . \nonumber
\end{equation}
Spiral $\Gamma$ is a wavy spiral if and only if $\nu\neq 0$. Analogous claim is valid for $x(t)=\widetilde{Y}_{\nu,\mu}(t)$, $t\in[\tau_0,\infty)$.
\end{theorem}

See Definition \ref{def_wavy_spiral} below for the wavy spiral. From Theorem \ref{tm_Bessel} directly follows our result about Bessel functions.

\begin{figure}[p]
\centering
\begin{subfigure}{.42\textwidth}
  \centering
  \includegraphics[width=.8\linewidth]{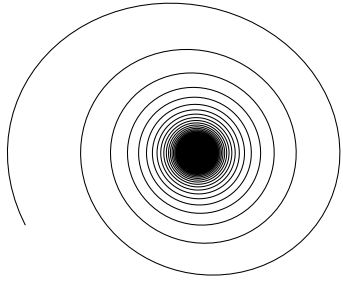}
  \caption{Curve $\Gamma_1$, $\dim_B\Gamma_1=\frac{20}{19}$}
\end{subfigure}%
\begin{subfigure}{.58\textwidth}
  \centering
  \includegraphics[width=.8\linewidth]{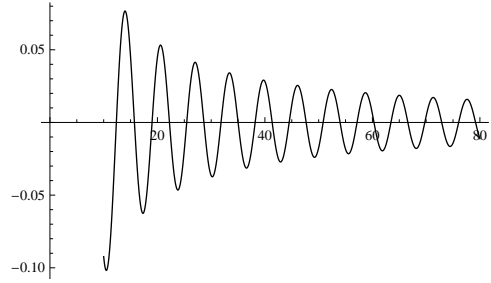}
  \caption{Graph of function $x_1(t)$}
\end{subfigure}
\caption{Curve $\Gamma_1=\{(x_1(t),\dot x_1(t)): t\in[10,\ty)\}$ and the graph of function $x_1(t)=\widetilde{J}_{5,0.2}(t)$. }
\label{fig:bessel_graf1}
\end{figure}

\begin{figure}[p]
\centering
\begin{subfigure}{.42\textwidth}
  \centering
  \includegraphics[width=.8\linewidth]{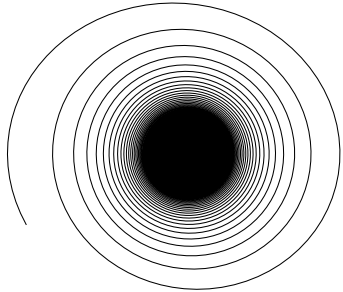}
  \caption{Curve $\Gamma_2$, $\dim_B\Gamma_2=\frac{4}{3}$}
\end{subfigure}%
\begin{subfigure}{.58\textwidth}
  \centering
  \includegraphics[width=.8\linewidth]{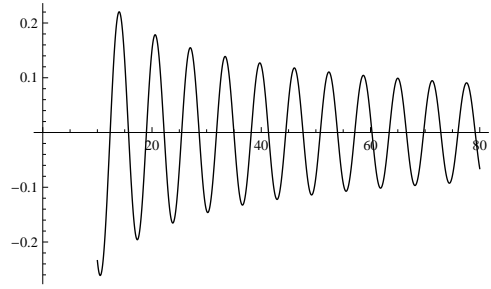}
  \caption{Graph of function $x_2(t)$}
\end{subfigure}
\caption{Curve $\Gamma_2=\{(x_2(t),\dot x_2(t)): t\in[10,\ty)\}$ and the graph of function $x_2(t)=\widetilde{J}_{5,1}(t)$. }
\label{fig:bessel_graf2}
\end{figure}

\begin{figure}[p]
\centering
\begin{subfigure}{.42\textwidth}
  \centering
  \includegraphics[width=.8\linewidth]{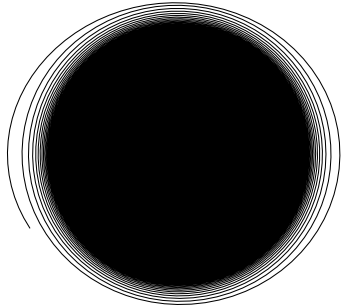}
  \caption{Curve $\Gamma_3$, $\dim_B\Gamma_3=\frac{20}{11}$}
\end{subfigure}%
\begin{subfigure}{.58\textwidth}
  \centering
  \includegraphics[width=.8\linewidth]{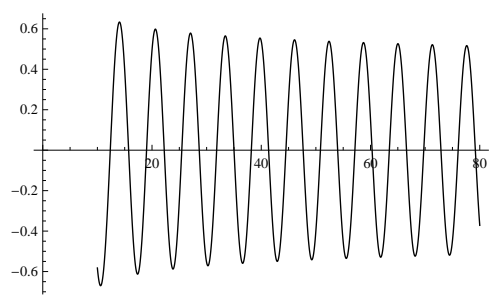}
  \caption{Graph of function $x_3(t)$}
\end{subfigure}
\caption{Curve $\Gamma_3=\{(x_3(t),\dot x_3(t)): t\in[10,\ty)\}$ and the graph of function $x_3(t)=\widetilde{J}_{5,1.8}(t)$. }
\label{fig:bessel_graf3}
\end{figure}

\begin{cor}{\label{cor_Bessel}} {\rm (Phase dimension of Bessel functions)} Phase dimension of Bessel functions is equal to $4/3$.
\end{cor}

{\it Proof.}
Follows directly from Theorem \ref{tm_Bessel} by taking $\mu=1$.
\qed

\bigskip

\begin{remark}

It is interesting to notice that value $4/3$ also appears as the box dimension of a trajectory near the origin of system
\begin{eqnarray}
\dot{r} & = & r(r^2+a_0),\quad a_0\in\eR \nonumber\\
\dot{\f} & = & 1 ,\nonumber
\end{eqnarray}
in polar coordinates, by taking $a_0=0$, see \cite{zuzu}. This system is the normal form of the Hopf bifurcation, which is a well known bifurcation of 1-parameter families of vector fields.

Another place where value $4/3$ appears is as the box dimension of the famous Euler spiral, or the Cornu spiral, also called the clothoid, defined parametricaly by
\begin{eqnarray}
x(t) & = & \int_0^t \cos(s^2)\,ds, \nonumber\\
y(t) & = & \int_0^t \sin(s^2)\,ds, \quad\textrm{where}\ t\in\eR , \nonumber
\end{eqnarray}
near $t=\pm\infty$, see \cite{clothoid}. The Euler spiral appears in various branches of mathematics, physics and engineering. Some applications are in optimal control theory, in the description of diffraction phenomena in optics and in robotics for pathfinding algorithms.

\end{remark}

In Theorem \ref{tm_Bessel} we have a spiral generated by generalized Bessel functions, see Figures \ref{fig:bessel_graf1}--\ref{fig:bessel_graf3}. To prove Theorem \ref{tm_Bessel} we need a new version of \cite[Theorem 5]{zuzu} with weaker assumptions on the spiral than in \cite{zuzu}.

\begin{theorem}{\label{novizuzu}} {\rm (Dimension of a piecewise smooth nonincreasing spiral)}
Let $f:[\f_1,\infty)\rightarrow(0,\infty)$ be a nonincreasing and radially decreasing function, also a continuous and piecewise continuously differentiable. We assume that the number of smooth pieces of $f$ in $[\f_1,\overline{\f}_1]$ is finite, for any $\overline{\f}_1>\f_1$.
Assume that there exist positive constants $\underline{m}$, $\overline{m}$, $\underline{a}$ and $M$ such that for all $\f\geq\f_1$,
$$
\underline{m}\f^{-\alpha}\leq f(\f)\leq \overline{m}\f^{-\alpha} ,
$$
$$
\underline{a}\f^{-\alpha-1}\leq f(\f)-f(\f+2\pi) ,
$$
and for all $\f$ where $f(\f)$ is differentiable,
$$
|f'(\f)|\leq M\f^{-\alpha-1} .
$$
Let $\Gamma$ be the graph of $r=f(\f)$ in polar coordinates. If $\alpha\in(0,1)$ then
$$
\dim_B\Gamma=\frac{2}{1+\alpha}.
$$
\end{theorem}

We also need the following Lemma~\ref{novalemazuzu} that is a generalization of \cite[Lemma~1]{zuzu} dealing with smooth spirals.

\begin{lemma}{\label{novalemazuzu}} {\rm (Excision property for piecewise smooth curves)}
Let $\Gamma$ be the image of a continuous and piecewise continuously differentiable
function $h:[\f_1,\infty)\to\mathbb{R}^2$ $($piecewise in the sense of Theorem \ref{novizuzu}$)$. Assume that $\underline{\dim}_B\Gamma>1$, $\Gamma_1:=h((\overline{\f}_1,\infty))$, for some fixed $\overline{\f}_1>\f_1$, and $h([\f_1,\overline{\f}_1])\bigcap\Gamma_1=\emptyset$. Then
$$
\underline{\dim}_B\Gamma_1=\underline{\dim}_B\Gamma,\quad \overline{\dim}_B\Gamma_1=\overline{\dim}_B\Gamma .
$$
\end{lemma}

Notice that $\Gamma$ from Lemma~\ref{novalemazuzu} is a piecewise smooth curve in $\eR^2$. Theorem~\ref{novizuzu} and Lemma~\ref{novalemazuzu} are introduced and proved in \cite{cswavy}.

\smallskip

Actually, we need a version of Theorem \ref{novizuzu} with even weaker assumptions on associated spiral. So, using Theorem \ref{novizuzu} and Lemma~\ref{novalemazuzu}, here we prove Theorem \ref{biliplema} that deals with a spiral $\Gamma'$ described by $r=f(\f)$, where $f$ is increasing on some parts, see Definitions \ref{def_wavy_function} and \ref{def_wavy_spiral}. We call this new property of $\Gamma'$ \emph{spiral  waviness}, which naturally arises in spirals generated by generalized Bessel functions in Theorem \ref{tm_Bessel}.

\begin{defin}\label{def_wavy_function} Let $r:[t_0,\infty)\rightarrow(0,\infty)$ be a $C^1$ function. Assume that $r'(t_0)\leq 0$. We say that $r=r(t)$ is a \emph{wavy function} if the sequence $(t_n)$ defined inductively by:
\begin{eqnarray}
t_{2k+1} & := & \inf\{t : t>t_{2k}, r'(t)>0\},\quad k\in\mathbb{N}_0 , \nonumber\\
t_{2k+2} & := & \inf\{t : t>t_{2k+1}, r(t)=r(t_{2k+1})\},\quad k\in\mathbb{N}_0 , \nonumber
\end{eqnarray}
is well-defined, and satisfies the \emph{waviness condition}:
\begin{equation}\label{def_SRC}
\left\{
\begin{minipage}{10cm}
\begin{itemize}
\item[(i)] The sequence $(t_n)$ is increasing and $t_n\to\infty$ as $n\to\infty$.
\item[(ii)] There exists $\varepsilon>0$, such that for all $k\in\mathbb{N}_0$ holds $t_{2k+1} - t_{2k} \geq\varepsilon$.
\item[(iii)] For all $k$ sufficiently large it holds
    $\mathop{\mathrm{osc}}\limits_{t\in[t_{2k+1},t_{2k+2}]} r(t) = o\left(t_{2k+1}^{-\alpha-1}\right)$, $\a\in(0,1)$,
\end{itemize}
\end{minipage}
\right.
\end{equation}
where $\mathop{\mathrm{osc}}\limits_{t\in I} r(t)=\max\limits_{t\in I} r(t) - \min\limits_{t\in I} r(t)$.
\end{defin}

Notice that $\min\limits_{t\in[t_{2k+1},t_{2k+2}]} r(t)=r(t_{2k+1})$. Condition $(i)$ means that the property of waviness of $r=r(\f)$ is global on the whole domain.
Condition~$(ii)$ is connected to an assumption of Lemma~\ref{tehnicka2}.
Condition $(iii)$ is a condition on a decay rate on the sequence of oscillations of $r$ on $I_k=[t_{2k+1},t_{2k+2}]$, for $k$ sufficiently large. Also, notice that condition $r'(t_0)\leq 0$ assures that $t_1$ is well-defined.

\begin{remark}\label{remark-tehnicki-uvjet}
Conditions (i) and (ii) in the waviness condition (\ref{def_SRC}) are not entirely independent. From (ii) and if $(t_n)$ is increasing follows that $t_n\to\infty$ as $n\to\infty$, but from (i) does not follow (ii). So, condition (ii) plus $(t_n)$ increasing is stronger than condition (i).
\end{remark}

\begin{defin}\label{def_wavy_spiral} Let a spiral $\Gamma'$, given in polar coordinates by $r=f(\f)$, where $f$ is a given function. If there exists increasing or decreasing function of class $C^1$, $\f=\f(t)$, such that $r(t)=f(\f(t))$ is a wavy function, then we say $\Gamma'$ is a \emph{wavy spiral}.
\end{defin}

Now, using Theorem \ref{novizuzu} and Lemma \ref{novalemazuzu} we prove the following Theorem \ref{biliplema}.

\begin{theorem}{\label{biliplema}} {\rm (Box dimension of a wavy spiral)}
Let $t_0>0$ and assume $r:[t_0,\infty)\rightarrow(0,\infty)$ is a wavy function. Assume that $\f:[t_0,\infty)\rightarrow[\f_0,\infty)$ is an increasing function of class $C^1$ such that $\f(t_0)=\f_0>0$ and there exists $\bar{\f}_0\in\eR$ such that
\begin{equation} \label{uvjetfikaot}
|(\f(t)-\bar{\f}_0)-(t-t_0)|\to 0\ \ \textrm{as}\ \ t\to\ty.
\end{equation}
Let $f:[\f_0,\infty)\rightarrow(0,\infty)$ be defined by $f(\f(t))=r(t)$. Assume that $\Gamma'$ is a spiral defined in polar coordinates by $r=f(\f)$, satisfying $(\ref{def_spiral})$.
Let $\a\in(0,1)$ is the same value as in $(\ref{def_SRC})${\rm(iii)} for wavy function $r$, and assume $\varepsilon'$ is such that $0<\varepsilon'<\varepsilon$, where $\varepsilon$ is defined by $(\ref{def_SRC})${\rm(ii)} for wavy function $r$. Assume that there exist positive constants $\underline{m}$, $\overline{m}$, $\underline{a}'$ and $M$ such that for all $\f\geq\f_0$,
\begin{equation}\label{biliplema_cond1}
\underline{m}\f^{-\alpha}\leq f(\f)\leq \overline{m}\f^{-\alpha} ,
\end{equation}
\begin{equation}\label{biliplema_cond3}
|f'(\f)|\leq M\f^{-\alpha-1} ,
\end{equation}
and for all $\triangle\f$, such that $\theta\leq\triangle\f\leq 2\pi+\theta$, there holds
\begin{equation}\label{biliplema_cond2b}
\underline{a}'\f^{-\alpha-1}\leq f(\f)-f(\f+\triangle\f) ,
\end{equation}
where $\theta:=\min\left\{\varepsilon',\pi\right\}$.

Then $\Gamma'$ is a wavy spiral and
$$
\dim_B\Gamma'=\frac{2}{1+\alpha}.
$$
\end{theorem}

{\it Proof.}

We will construct a new spiral $\widetilde{\Gamma}$ that is nonincreasing and close to $\Gamma'$ in the sense that there exists bi-Lipschitz map $F:\Gamma'\to\widetilde{\Gamma}$ (ie, there exist constants $\overline{K}_1\in(0,1)$ and $\overline{K}_2>1$ such that
$$
\overline{K}_1\mathrm{d}(X',Y')\leq \mathrm{d}(F(X'),F(Y'))\leq \overline{K}_2\mathrm{d}(X',Y') ,
$$
for every $X',Y'\in\Gamma'$, where $\mathrm{d}(X',Y')$ is Euclidian distance between points $X'$ and $Y'$). Also, we will construct this map $F$. We know from \cite[p.\ 44]{falc} that then $\dim_B\Gamma'=\dim_B\widetilde{\Gamma}$. Because of Lemma \ref{novalemazuzu} we have to construct $\widetilde{\Gamma}$ and establish bi-Lipschitz map $F$ only for $\f_0$ and equivalently $t_0$ sufficiently large. This means we can increase value $t_0$ to $t_{2k_0}$ and take $\f_0=\f(t_{2k_0})$ for any $k_0\in\eN$, where sequence $(t_n)$ is from Definition \ref{def_wavy_function}.

Let $\widetilde{\Gamma}$ be defined by $r=\widetilde{f}(\f)$, in polar coordinates. Informally, $\widetilde{\Gamma}$ is composed of alternating radially nonincreasing parts of spiral $\Gamma'$ (with respect to parameter $t$) and circular arcs that substitute the parts where $\Gamma'$ radially increases and then decreases until it intersects the associated arc. We see that $\widetilde{\Gamma}$ is radially nonincreasing.

Formally, we define $\widetilde{\Gamma}$ by parts as follows:
\begin{eqnarray}
\widetilde{\Gamma}|_{[t_{2k},t_{2k+1}]} & := & \Gamma'|_{[t_{2k},t_{2k+1}]},\quad\textrm{for every}\ k\in\eN_0, \nonumber\\
\widetilde{\Gamma}|_{[t_{2k+1},t_{2k+2}]} & := & \left\{\left(\widetilde{r}(t),\f(t)\right) : \widetilde{r}(t):=r(t_{2k+1}),\ t\in[t_{2k+1},t_{2k+2}] \right\}, \nonumber\\
& & \textrm{for every}\ k\in\eN_0, \nonumber
\end{eqnarray}
in polar coordinates.

As $r$ is a wavy function (see Definition \ref{def_wavy_function}) and $\Gamma'$ is a spiral, satisfying $(\ref{def_spiral})$, we see that $\widetilde{\Gamma}$ is well defined and is also a spiral, satisfying $(\ref{def_spiral})$. Notice that $\widetilde{f}\left(\f(t)\right)=\widetilde{r}(t)$ and $\widetilde{r}(t)\leq r(t)$, for all $t\in[t_0,\infty)$.

Assume $T'\in\Gamma'$, which means $T'$ is given by $(r(t_{T'}),\f(t_{T'}))$ for some value $t_{T'}$, in polar coordinates. Now we define
\begin{equation}
F(T'):=\left(\widetilde{r}(t_{T'}),\f(t_{T'})\right) , \nonumber
\end{equation}
in polar coordinates. We see that map $F$ radially maps points from spiral $\Gamma'$ to spiral $\widetilde{\Gamma}$.

We will use $f(T')$ and $\widetilde{f}(T)$ to denote radius, in polar coordinates, of points $T'\in\Gamma'$ and $T\in\widetilde{\Gamma}$, respectively. With $\f_T$ we denote argument, in polar coordinates, corresponding to points $T'\in\Gamma'$ or $T=F(T')\in\widetilde{\Gamma}$. With $O$ we denote the origin of the phase plane.

Now, we prove that map $F$ is a bi-Lipschitz map. Notice this is equivalent to the claim that $F^{-1}:\widetilde{\Gamma}\to\Gamma'$ is a bi-Lipschitz map.

We have five cases depending on the relative position of points $X,Y\in\widetilde{\Gamma}$. For each case we will prove that $F^{-1}$ is a bi-Lipschitz map.

\smallskip

\textbf{Case A.} Assume point $X\in\widetilde{\Gamma}|_{[t_{2k},t_{2k+1}]}$, for some $k\in\mathbb{N}_0$, and point $Y\in\widetilde{\Gamma}|_{[t_{2l},t_{2l+1}]}$, for some $l\in\mathbb{N}_0$. Map $F^{-1}$ is identity map on $\widetilde{\Gamma}|_{[t_{2k},t_{2k+1}]}\bigcup \widetilde{\Gamma}|_{[t_{2l},t_{2l+1}]}$ so $F^{-1}$ is, trivially, a bi-Lipschitz map on $\widetilde{\Gamma}|_{[t_{2k},t_{2k+1}]}\bigcup \widetilde{\Gamma}|_{[t_{2l},t_{2l+1}]}$, for every $k, l\in\mathbb{N}_0$.

\smallskip

\begin{figure}
\begin{center}
\includegraphics{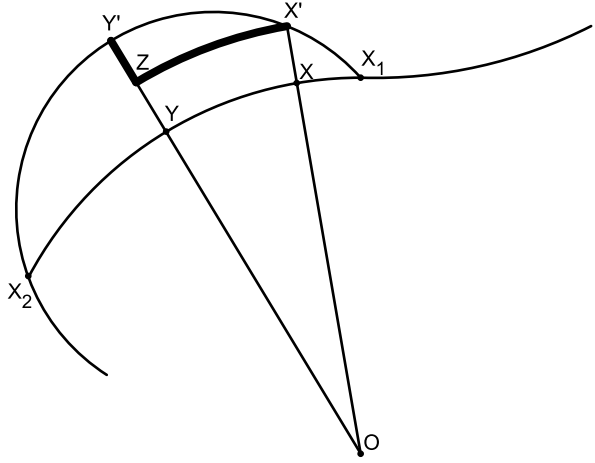}
\end{center}
\caption{Geometry in the proof of Theorem \ref{biliplema}, \bf{Case B}.}
\label{slika-caseB}
\end{figure}

\textbf{Case B.} Assume points $X, Y \in\widetilde{\Gamma}|_{[t_{2k+1},t_{2k+2}]}$, for some $k\in\mathbb{N}_0$, and $|\f_X-\f_Y|<\theta$. Without loss of generality we can assume that $\f_X<\f_Y$. Let us denote with $X'=F^{-1}(X)$, $Y'=F^{-1}(Y)$, and with $X_1$ ($X_2$) the point on spiral $\widetilde{\Gamma}$ corresponding to parameter $t_{2k+1}$ ($t_{2k+2}$) and angle $\f_{2k+1}:=\f(t_{2k+1})$ ($\f_{2k+2}:=\f(t_{2k+2})$), see Figure \ref{slika-caseB}. As $\f$ is increasing, we see that $\f_{2k+1}<\f_{2k+2}$.

Using elementary geometry regarding isosceles triangle $XYO$ and triangle $X'Y'O$ (see Figure \ref{slika-caseB}) we conclude that
$$d(X,Y) \leq d(X',Y') = d(F^{-1}(X),F^{-1}(Y)) .$$

Now, let us assume that $f(\f_Y)\geq f(\f_X)$ and let point $Z$ on the line $\overline{OY'}$ be such that $d(O,Z)=d(O,X')$. We see that
\begin{eqnarray}
d(F^{-1}(X),F^{-1}(Y))  & =    & d(X',Y') \leq d(X',Z)+d(Z,Y') \nonumber\\
                        & \leq & |\f_X-\f_Y|\cdot d(O,X') + |f(X')-f(Y')| \nonumber\\
                        & =    & |\f_X-\f_Y|\cdot f(\f_X) + |f(\f_X)-f(\f_Y)| \nonumber\\
                        & =    & |\f_X-\f_Y|\cdot f(\f_X) + |f'(\xi)||\f_X-\f_Y|, \nonumber\\
                        &      & \quad \xi\in[\f_X,\f_Y]\subseteq[\f_{2k+1},\f_{2k+2}] \nonumber\\
                        & \leq & |\f_X-\f_Y|\left[f(\f_X)+\sup\limits_{\xi\in[\f_{2k+1},\f_{2k+2}]}|f'(\xi)|\right] \nonumber\\
                        & =    & |\f_X-\f_Y|f(\f_{2k+1})\cdot\frac{f(\f_X)+\sup\limits_{\xi\in[\f_{2k+1},\f_{2k+2}]}|f'(\xi)|}{f(\f_{2k+1})} . \nonumber
\end{eqnarray}

Analogously, if $f(\f_Y)< f(\f_X)$ we conclude that
$$
d(F^{-1}(X),F^{-1}(Y)) \leq |\f_X-\f_Y|f(\f_{2k+1})\cdot\frac{f(\f_Y)+\sup\limits_{\xi\in[\f_{2k+1},\f_{2k+2}]}|f'(\xi)|}{f(\f_{2k+1})} .
$$

As $\theta\leq\pi$, using Proposition \ref{geometrijska} we see that
$$
|\f_X-\f_Y|f(\f_{2k+1})\leq \frac{\pi}{2}d(X,Y) .
$$

As function $f$ satisfies conditions (\ref{biliplema_cond1}) and (\ref{biliplema_cond3}), and as we can take $\f_0$ sufficiently large such that $M\f_{2k+1}^{-\alpha-1}\leq\overline{m}|\f_{2k+1}|^{-\alpha}$, we compute
\begin{eqnarray}
\frac{f(\f_X)+\sup\limits_{\xi\in[\f_{2k+1},\f_{2k+2}]}|f'(\xi)|}{f(\f_{2k+1})} & \leq & \frac{\overline{m}\f_X^{-\alpha}+\sup\limits_{\xi\in[\f_{2k+1},\f_{2k+2}]}(M\xi^{-\alpha-1})}{\underline{m}\f_{2k+1}^{-\alpha}} \nonumber\\
& \leq & \frac{\overline{m}\f_{2k+1}^{-\alpha}+M\f_{2k+1}^{-\alpha-1}}{\underline{m}\f_{2k+1}^{-\alpha}} \leq \frac{\overline{m}\f_{2k+1}^{-\alpha}+\overline{m}\f_{2k+1}^{-\alpha}}{\underline{m}\f_{2k+1}^{-\alpha}} \nonumber\\
& \leq & 2\frac{\overline{m}}{\underline{m}} .\nonumber
\end{eqnarray}
Analogously, we get
$$
\frac{f(\f_Y)+\sup\limits_{\xi\in[\f_{2k+1},\f_{2k+2}]}|f'(\xi)|}{f(\f_{2k+1})} \leq 2\frac{\overline{m}}{\underline{m}} .
$$

Finally, we conclude that
$$
d(F^{-1}(X),F^{-1}(Y)) \leq \frac{\pi}{2}d(X,Y)\cdot 2\frac{\overline{m}}{\underline{m}} = \pi\frac{\overline{m}}{\underline{m}}d(X,Y) ,
$$
so $F^{-1}$ is a bi-Lipschitz map on $\widetilde{\Gamma}|_{[t_{2k+1},t_{2k+2}]}$, for every $k\in\mathbb{N}_0$.

\smallskip

\textbf{Case C1.} Assume point $X \in\widetilde{\Gamma}|_{[t_{2k+2},t_{2k+3}]}$ and point $Y \in\widetilde{\Gamma}|_{[t_{2k+1},t_{2k+2}]}$, for some $k\in\mathbb{N}_0$, and $|\f_X-\f_Y|<\theta$. Let us denote with $Z \in\widetilde{\Gamma}$ the point corresponding to parameter $t_{2k+2}$ and angle $\f_Z$. Notice that $\f_X\geq \f_Z\geq \f_Y$.

\begin{figure}
\begin{center}
\includegraphics[width=8cm]{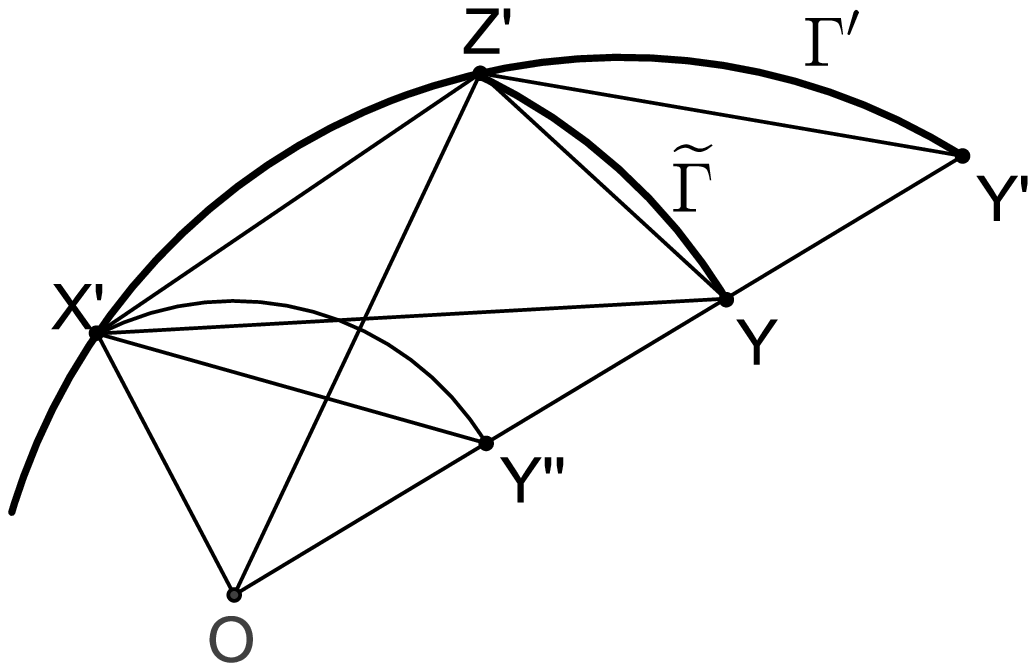}
\end{center}
\caption{Geometry in the proof of Theorem \ref{biliplema}, \bf{Case C1}.}
\label{slika-caseC1}
\end{figure}

Let us denote with $X'=F^{-1}(X)$, $Y'=F^{-1}(Y)$ and $Z'=F^{-1}(Z)$, see Figure \ref{slika-caseC1}. Notice that $X=X'$, $Z=Z'$, which gives $d(X,Z)=d(X',Z')$. Now,
$$
\widetilde{f}(X)=d(O,X)=d(O,X')=f(X') ,
$$
$$
\widetilde{f}(Z)=d(O,Z)=d(O,Z')=f(Z') ,
$$
$$
d(O,Y)=\widetilde{f}(Y)\leq f(Y')=d(O,Y') ,
$$
and notice that
$$
d(O,Y)=\widetilde{f}(Y)=\widetilde{f}(Z)=d(O,Z) ,
$$
$$
f(Z')=\widetilde{f}(Z)=\widetilde{f}(Y)\leq f(Y') .
$$
From
$$
d(O,X)=\widetilde{f}(X)\leq \widetilde{f}(Z)=\widetilde{f}(Y)=d(0,Y)
$$
easily follows (see Figure \ref{slika-caseC1}) that the angle between line $\overline{XY}$ and line $\overline{YY'}$ in triangle $XYY'$ is greater than or equal to $\pi/2$, so using elementary geometry we see that
$$
d(X,Y)\leq d(X',Y')=d(F^{-1}(X),F^{-1}(Y)) .
$$
On the other hand, from \textbf{Case B}, it follows $d(Z',Y')\leq \pi\frac{\overline{m}}{\underline{m}} d(Z,Y)$.

Now, we compute
\begin{eqnarray}
&    & d(F^{-1}(X),F^{-1}(Y)) = d(X',Y') \leq d(X',Z')+d(Z',Y') \nonumber\\
&\leq& d(X,Z)+\pi\frac{\overline{m}}{\underline{m}} d(Z,Y)\leq \pi\frac{\overline{m}}{\underline{m}}\left[d(X,Z)+d(Z,Y)\right] \nonumber\\
&\leq& \pi\frac{\overline{m}}{\underline{m}}\left[|\f_X-\f_Z|\max\{d(O,X),d(O,Z)\}+(d(O,Z)-d(O,X))\right. \nonumber\\
&+   & \left.|\f_Z-\f_Y|\max\{d(O,Z),d(O,Y)\}\right] \nonumber\\
&=   & \pi\frac{\overline{m}}{\underline{m}}\left[(|\f_X-\f_Z|+|\f_Z-\f_Y|)\widetilde{f}(Z)+(f(Z')-f(X'))\right] \nonumber\\
&=   & \pi\frac{\overline{m}}{\underline{m}}\left[|\f_X-\f_Y|f(Z')+(f(\f_Z)-f(\f_X))\right] \nonumber\\
&\leq& \pi\frac{\overline{m}}{\underline{m}}\left[|\f_X-\f_Y|f(\f_Z)+|\f_X-\f_Z|\sup\limits_{\xi\in[\f_Z,\f_X]}|f'(\xi)|\right] \nonumber\\
&\leq& \pi\frac{\overline{m}}{\underline{m}}|\f_X-\f_Y|\left[f(\f_Y)+\sup\limits_{\xi\in[\f_Y,\f_X]}|f'(\xi)|\right] \nonumber\\
&=   & \pi\frac{\overline{m}}{\underline{m}} \cdot |\f_X-\f_Y|f(\f_X) \cdot \frac{f(\f_Y)}{f(\f_X)} \cdot \frac{f(\f_Y)+\sup\limits_{\xi\in[\f_Y,\f_X]}|f'(\xi)|}{f(\f_Y)} \nonumber
\end{eqnarray}

Let point $Y''\in\overline{OY}$ be such that $d(O,Y'')=d(O,X)$. Notice that $d(O,Y'')\leq d(O,Y)$. As before, using elementary geometry regarding isosceles triangle $XOY''$ and triangle $XOY$ (see Figure \ref{slika-caseC1}) we conclude that $d(X,Y'')\leq d(X,Y)$. Using Proposition \ref{geometrijska} on isosceles triangle $XOY''$, we get that
$$
|\f_X-\f_Y|f(\f_X)\leq \frac{\pi}{2}d(X,Y'')\leq  \frac{\pi}{2}d(X,Y) .
$$

We can take $\f_0$ sufficiently large such that, using $\f_X\geq \f_Y\geq \f_0$, we get
\begin{eqnarray}
\frac{f(\f_Y)}{f(\f_X)} & \leq & \frac{\overline{m}\f_Y^{-\alpha}}{\underline{m}\f_X^{-\alpha}} = \frac{\overline{m}\f_X^{\alpha}}{\underline{m}\f_Y^{\alpha}} \leq \frac{\overline{m}(\f_Y+|\f_X-\f_Y|)^{\alpha}}{\underline{m}\f_Y^{\alpha}} \nonumber\\
& \leq & \frac{\overline{m}(\f_Y+\theta)^{\alpha}}{\underline{m}\f_Y^{\alpha}} = \frac{\overline{m}\f_Y^{\alpha}(1+O(\f_Y^{-1}))}{\underline{m}\f_Y^{\alpha}} \leq 2\frac{\overline{m}}{\underline{m}} . \nonumber
\end{eqnarray}

Analogously as in \textbf{Case B}, we can take $\f_0$ sufficiently large such that
$$
\frac{f(\f_Y)+\sup\limits_{\xi\in[\f_Y,\f_X]}|f'(\xi)|}{f(\f_Y)} \leq 2\frac{\overline{m}}{\underline{m}} .
$$

Finally, we see that
$$
d(F^{-1}(X),F^{-1}(Y)) \leq \pi\frac{\overline{m}}{\underline{m}} \cdot \frac{\pi}{2}d(X,Y) \cdot 2\frac{\overline{m}}{\underline{m}} \cdot 2\frac{\overline{m}}{\underline{m}} = 2\pi^2 \left(\frac{\overline{m}}{\underline{m}}\right)^3 d(X,Y) .
$$

\smallskip

\textbf{Case C2.} Assume point $X \in\widetilde{\Gamma}|_{[t_{2k+1},t_{2k+2}]}$ and point $Y \in\widetilde{\Gamma}|_{[t_{2k},t_{2k+1}]}$, for some $k\in\mathbb{N}_0$, and $|\f_X-\f_Y|<\theta$. Proof is analogous to proof of \textbf{Case C1}.

\smallskip

\textbf{Conclusion of Cases A--C2.} Using (\ref{def_SRC})(ii) and (\ref{uvjetfikaot}), as $t_{2k+1}-t_{2k}\geq\varepsilon>\theta$ for all $k\in\eN_0$, it follows that, for $t_0$ and equivalently $\f_0$ sufficiently large, cases \textbf{A} to \textbf{C2} cover all possibilities for relative positions of points $X,Y\in\widetilde{\Gamma}$ such that $|\f_X-\f_Y|<\theta$. For $\f_0$ sufficiently large, it follows that $F^{-1}$ is a bi-Lipschitz map, for all points $X,Y\in\widetilde{\Gamma}$ such that $|\f_X-\f_Y|<\theta$.

\smallskip

\textbf{Case D.} Assume points $X, Y \in\widetilde{\Gamma}$ and $|\f_X-\f_Y|\geq\theta$. Without loss of generality we can assume that $\f_X<\f_Y$. Let us denote with $X'=F^{-1}(X)$ and $Y'=F^{-1}(Y)$.

First, we consider the case when $|\f_X-\f_Y| < 2\pi+\theta$. From (\ref{biliplema_cond2b}) we compute
\begin{eqnarray} \label{ddonjaograda}
d(X',Y') & \geq & d(X',O) - d(Y',O) = f(\f_X) - f(\f_Y) \\
& = & f(\f_X) - f(\f_X+(\f_Y-\f_X)) \geq \underline{a}'\f_X^{-\alpha-1} . \nonumber
\end{eqnarray}

\setlength{\unitlength}{.5\linewidth}
\begin{figure}
\begin{center}
\begin{picture}(1,1.15)
\put(0,0){\includegraphics[width=.5\linewidth]{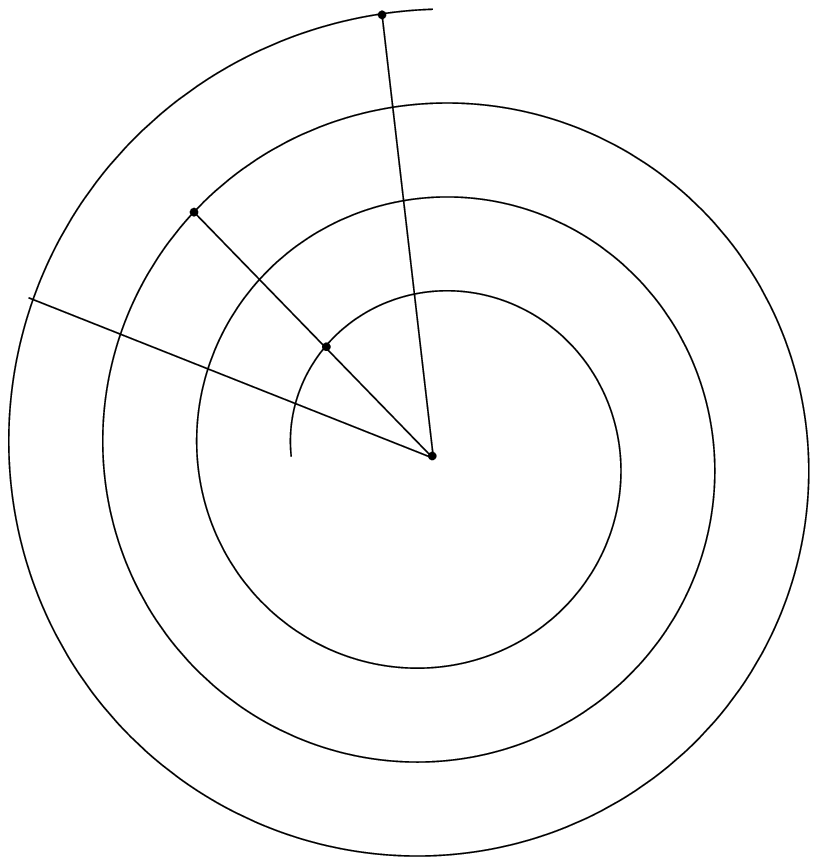}}
\put(0.45,1.12){$X'$}
\put(0.375,0.745){$Y'$}
\put(0.215,0.905){$Y''$}
\put(0.5,0.52){$O$}
\end{picture}
\end{center}
\caption{Geometry in the proof of Theorem \ref{biliplema}, \bf{Case D}.}
\label{slika-caseD}
\end{figure}

Next, we consider the case when $|\f_X-\f_Y| \geq 2\pi+\theta$. We extend line $\overline{OY'}$ such that it intersects spiral $\Gamma'$ at point $Y''\in\Gamma'$ such that $\f_{Y''}\in[\theta,2\pi+\theta)$ (see Figure \ref{slika-caseD}). It is easy to see that point $Y''$ is well defined. We see that $f(Y')<f(Y'')$, because $f$ is radially decreasing, see (\ref{def_spiral}).

Using (\ref{ddonjaograda}), we compute
\begin{eqnarray} \label{ddonjaograda2}
d(X',Y') & \geq & d(X',O) - d(Y',O) = f(X') - f(Y') \\
& \geq & f(X') - f(Y'') = f(\f_X) - f(\f_{Y''}) \nonumber\\
& \geq &\underline{a}'\f_X^{-\alpha-1} . \nonumber
\end{eqnarray}

From triangle inequality we know that
\begin{eqnarray} \label{nejtrokuta}
d(X',Y') & \leq & d(X',X)+d(X,Y)+d(Y,Y') , \\
d(X,Y) & \leq & d(X,X')+d(X',Y')+d(Y',Y) . \nonumber
\end{eqnarray}
Now, we define $\f_n:=\f(t_n)$, $\forall n\in\mathbb{N}_0$. Notice that from (\ref{uvjetfikaot}) follows $\f_n \simeq t_n$ as $n\to\infty$.

It is easy to see that $t_{2k+2}-t_{2k+1}<2\pi$ for every $k\in\eN_0$. (On the contrary, $\widetilde{\Gamma}$ would not be a spiral, see (\ref{def_spiral}), because it would have a self intersection.) Using (\ref{uvjetfikaot}) we see that $\f_{2k+2}-\f_{2k+1}<2\pi+1$ for every $k\in\eN_0$. (For the upper bound we could take any number larger than $2\pi$.)

Now, assume $Z'\in\Gamma'$ and $Z:=F(Z')$. There exists $k_Z\in\eN_0$ such that $\f_Z\in[\f_{2k_Z},\f_{2k_Z+2}]$. If $\f_Z\in[\f_{2k_Z},\f_{2k_Z+1}]$ then $Z\in\widetilde{\Gamma}|_{[t_{2k_Z},t_{2k_Z+1}]}$, so
$$
d(Z,Z')=0 .
$$
If $\f_Z\in[\f_{2k_Z+1},\f_{2k_Z+2}]$ then $Z\in\widetilde{\Gamma}|_{[t_{2k_Z+1},t_{2k_Z+2}]}$, so (\ref{def_SRC})(iii) gives
\begin{eqnarray}
d(Z,Z') & \leq & \sup\limits_{t\in[t_{2k_Z+1},t_{2k_Z+2}]}\{r(t)-r(t_{2k_Z+1})\} \nonumber\\
& = & \mathop{\mathrm{osc}}\limits_{t\in[t_{2k_Z+1},t_{2k_Z+2}]} r(t) = o\left(t_{2k_Z+1}^{-\alpha-1}\right) = o\left(\f_{2k_Z+1}^{-\alpha-1}\right) = o\left(\f_Z^{-\alpha-1}\right), \nonumber
\end{eqnarray}
for $\f_0$ sufficiently large. So, for every $Z'\in\Gamma'$ and $Z:=F(Z')$ is
\begin{equation}\label{dgornjaograda2}
0\leq d(Z,Z')=o\left(\f_Z^{-\alpha-1}\right) ,
\end{equation}
for $\f_0$ sufficiently large. From (\ref{nejtrokuta}) it follows
$$
d(X',Y')-d(X,X')-d(Y,Y') \leq d(X,Y) \leq d(X',Y')+d(X,X')+d(Y,Y') .
$$
Using (\ref{dgornjaograda2}) we get
\begin{eqnarray}
& & d(X',Y')-|o\left(\f_X^{-\alpha-1}\right)|-|o\left(\f_Y^{-\alpha-1}\right)| \leq d(X,Y) , \nonumber\\
& & d(X,Y) \leq d(X',Y')+|o\left(\f_X^{-\alpha-1}\right)|+|o\left(\f_Y^{-\alpha-1}\right)| , \nonumber
\end{eqnarray}
for $\f_0$ sufficiently large, which is equivalent to
$$
\left(1-\frac{|o\left(\f_X^{-\alpha-1}\right)|}{d(X',Y')}\right)d(X',Y') \leq d(X,Y) \leq \left(1+\frac{|o\left(\f_X^{-\alpha-1}\right)|}{d(X',Y')}\right)d(X',Y') ,
$$
for $\f_0$ sufficiently large. (Notice that as $\f_X<\f_Y$, by abusing notation we have $o\left(\f_Y^{-\alpha-1}\right)\subseteq o\left(\f_X^{-\alpha-1}\right)$, for $\f_0$ sufficiently large.)  Respecting the lower bound on $d(X',Y')$, (\ref{ddonjaograda}) and (\ref{ddonjaograda2}), we get
$$
(1-\frac{|o(1)|}{\underline{a}'})d(X',Y') \leq d(X,Y) \leq (1+\frac{|o(1)|}{\underline{a}'})d(X',Y') ,
$$
for $\f_0$ sufficiently large, which means there exist constants $K_3\in(0,1)$ and $K_4>1$ such that
$$
K_3 d(X',Y')\leq d(X,Y) \leq K_4 d(X',Y') ,
$$
for $\f_0$ sufficiently large. So finally, for $\f_0$ sufficiently large it follows that $F^{-1}$ is a bi-Lipschitz map on $\widetilde{\Gamma}$ for points $X,Y\in\widetilde{\Gamma}$ such that $|\f_X-\f_Y|\geq\theta$.

\smallskip

From cases \textbf{A} to \textbf{D} it follows that $F^{-1}$ is a bi-Lipschitz map on spiral $\widetilde{\Gamma}$, for $\f_0$ sufficiently large. On the other hand, using (\ref{def_SRC})(ii) we see that the number of smooth pieces of function $\widetilde{f}(\f)$ on any bounded interval is finite and using $\f_{2k+2}-\f_{2k+1}<2\pi+1$, for every $k\in\eN_0$, we see that spiral $\widetilde{\Gamma}$, for $\f_0$ sufficiently large, satisfies all assumptions of Theorem \ref{novizuzu}. Using Lemma \ref{novalemazuzu} and as $F$ is a bi-Lipschitz map on $\Gamma'$, for $\f_0$ sufficiently large, we finally get
$$
\dim_B\Gamma'=\dim_B\widetilde{\Gamma} = \frac{2}{1+\alpha} .
$$

\qed

\smallskip

Finally, using Theorem \ref{biliplema} and Lemma \ref{novalemazuzu} we provide proof of Theorem \ref{tm_Bessel} about the box dimension of a spiral generated by generalized Bessel functions. This proof consists of checking out the conditions of Theorem \ref{biliplema}. The following lemma makes this checking easier.

\begin{lemma}{\label{tehnicka2}}
Let $\a\in(0,1)$ and
$$
r(t)=p(t)\sqrt{1 + \frac{p'(t)}{p(t)}f(t) + g(t)},\quad t\in[t_0,\ty),\ t_0>0,
$$
where $p(t)\sim_1t^{-\a}$ as $t\to\ty$, $|f(t)|\leq 1$ for all $t\in[t_0,\ty)$ and $g(t)\in O\left(t^{-2}\right)$ as $t\to\infty$. Let $C\in\eR$ and assume that $t(\f)=\f+C+O(\f^{-1})$ as $\f\to\ty$. Let $\triangle\f>1$.

Then there exists constant $k>0$, independent of $\f$ and $\triangle\f$, such that for all $\f$ sufficiently large it holds
\begin{equation} \label{prop_tehn2_uvjet}
r(t(\f))-r(t(\f+\triangle\f)) \geq k\f^{-\a-1}(1+O(\f^{-1})) .
\end{equation}

\end{lemma}
The proof of Lemma \ref{tehnicka2} is omitted, because it is very similar to the proof of \cite[Lemma 3]{cswavy}.

\begin{remark} \label{remark_prop_tehn2}
To the contrary to constant $k$, the constant hidden in term $O\left(\f^{-1}\right)$ in (\ref{prop_tehn2_uvjet}) is not independent of $\triangle\f$. This is the reason we must prescribe upper bound on term $\triangle\f$, in parts of the proof, where we will use Lemma~\ref{tehnicka2}.
\end{remark}

\medskip

Throughout the proof of Theorem \ref{tm_Bessel} we will use Hankel's asymptotic expansions of Bessel functions $J_{\nu}(t)$ and $Y_{\nu}(t)$ for large $t$, see \cite[p.~120]{lebedev}, as follows:
\begin{eqnarray}
J_{\nu}(t) & = & \left(\frac{2}{\pi t}\right)^{\frac{1}{2}}\left[P_{\nu}(t)\cos\chi-Q_{\nu}(t)\sin\chi\right] , \label{bessel_J} \\
Y_{\nu}(t) & = & \left(\frac{2}{\pi t}\right)^{\frac{1}{2}}\left[P_{\nu}(t)\sin\chi+Q_{\nu}(t)\cos\chi\right] , \label{bessel_Y}
\end{eqnarray}
where $\nu\in\eR$ and $\chi=t-\left(\frac{1}{2}\nu+\frac{1}{4}\right)\pi$. $P_{\nu}(t)$ and $Q_{\nu}(t)$ are given by
\begin{eqnarray}
P_{\nu}(t) & = & \sum\limits_{k=0}^{N} (-1)^k \frac{(\nu,2k)}{(2t)^{2k}} + O\left(t^{-2N-2}\right), \label{fun_P}\\
Q_{\nu}(t) & = & \sum\limits_{k=0}^{N} (-1)^k \frac{(\nu,2k+1)}{(2t)^{2k+1}} + O\left(t^{-2N-3}\right) \label{fun_Q},
\end{eqnarray}
as $t\to\infty$, which are asymptotic expansions to $N$ terms, where we introduce the notation
\begin{eqnarray}
(\nu,k) & = & \frac{(-1)^k}{k!}(\frac{1}{2}-\nu)_k(\frac{1}{2}+\nu)_k \nonumber\\
& = & \frac{(4\nu^2-1)(4\nu^2-3^2)\cdots(4\nu^2-(2k-1)^2)}{2^{2k}k!} , \nonumber\\
(\nu,0) & = & 1 . \nonumber
\end{eqnarray}
We will also use asymptotic expansions of derivatives
\begin{eqnarray}
\frac{d}{dt}J_{\nu}(t) & = & \left(\frac{2}{\pi t}\right)^{\frac{1}{2}} \left[-R_{\nu}(t)\sin\chi - S_{\nu}(t)\cos\chi\right] , \label{der_bessel_J} \\
\frac{d}{dt}Y_{\nu}(t) & = & \left(\frac{2}{\pi t}\right)^{\frac{1}{2}} \left[R_{\nu}(t)\cos\chi - S_{\nu}(t)\sin\chi\right] , \label{der_bessel_Y}
\end{eqnarray}
where $R_{\nu}(t)$ and $S_{\nu}(t)$ are given by
\begin{eqnarray}
R_{\nu}(t) & = & \sum\limits_{k=0}^{N} (-1)^k \frac{4\nu^2+16k^2-1}{4\nu^2-(4k-1)^2}\,\frac{(\nu,2k)}{(2t)^{2k}} + O\left(t^{-2N-2}\right), \label{fun_R}\\
S_{\nu}(t) & = & \sum\limits_{k=0}^{N} (-1)^k \frac{4\nu^2+4(2k+1)^2-1}{4\nu^2-(4k+1)^2}\,\frac{(\nu,2k+1)}{(2t)^{2k+1}} + O\left(t^{-2N-3}\right) \label{fun_S} ,
\end{eqnarray}
as $t\to\infty$.

Notice, that because $J_{\nu}(t)$ and $Y_{\nu}(t)$ are analytic functions for $t>0$, see \cite{lebedev}, it is easily verified that $P_{\nu}(t)$, $Q_{\nu}(t)$, $R_{\nu}(t)$ and $S_{\nu}(t)$ are analytic functions for $t>0$. This is a sufficient condition for asymptotic power series of these functions to be differentiable term by term, see \cite[1.6]{erdelyi}. This property will enable us to easily compute expressions containing higher derivatives of those functions.

\smallskip

{\it Proof of Theorem \ref{tm_Bessel}.}

\smallskip

First, we assume that $x(t):=\widetilde{J}_{\nu,\mu}(t)$, for every $t\in[\tau_0,\infty)$. The proof is presented in a series of steps:

\smallskip

\emph{Step 1.} (The box dimension is invariant with respect to mirroring of a spiral.) We will prove the equivalent claim, that planar curve $\Gamma'=\{(x(t),-\dot x(t)) : t\in[\tau_0,\ty)\}$ is a spiral, see (\ref{def_spiral}), defined by $r=f(\f)$, $\f\in[\phi_0, \ty)$, near the origin, satisfying $f(\f)\simeq\f^{-\frac{2-\mu}{2}}$, in polar coordinates, near the origin, and $\dim_B\Gamma'=\frac{4}{4-\mu}$. It is easy to see that curve $\Gamma$ is a mirror image of curve $\Gamma'$, with respect to the $x$-axis, hence $\Gamma$ is a spiral. Notice that if curve $\Gamma'$ is a wavy spiral then curve $\Gamma$ is also a wavy spiral. Reflecting with respect to the $x$-axis in the plane is an isometric map. As the isometric map is bi-Lipschitzian and therefore it preserves box dimensions, see \cite[p. 44]{falc}, we see that $\dim_B\Gamma = \dim_B\Gamma'=\frac{4}{4-\mu}$.

\smallskip

\emph{Step 2.} (Checking condition (\ref{uvjetfikaot}).)
From
\begin{eqnarray} \label{def_fazni_sustav}
 x(t) & = & \widetilde{J}_{\nu,\mu}(t) = t^{\frac{\mu-1}{2}} J_{\widetilde{\nu}}(t) , \nonumber  \\
 \dot x (t)&=&  \frac{d}{dt}\widetilde{J}_{\nu,\mu}(t) = \frac{\mu-1}{2}t^{\frac{\mu-3}{2}} J_{\widetilde{\nu}}(t) + t^{\frac{\mu-1}{2}} \frac{d}{dt}J_{\widetilde{\nu}}(t) ,
\end{eqnarray}
where $\widetilde{\nu}=\sqrt{\left(\frac{\mu-1}{2}\right)^2+\nu^2}$, using (\ref{bessel_J}--\ref{fun_S}), we compute
\bgeq\label{fip}
\tan \f(t)=-\frac{\dot x(t)}{x(t)}=\frac{\left(R_{\widetilde{\nu}}(t)+\frac{\mu-1}{2t}Q_{\widetilde{\nu}}(t)\right)\sin\widetilde{\chi} + \left(S_{\widetilde{\nu}}(t)-\frac{\mu-1}{2t}P_{\widetilde{\nu}}(t)\right)\cos\widetilde{\chi}} {P_{\widetilde{\nu}}(t)\cos\widetilde{\chi}-Q_{\widetilde{\nu}}(t)\sin\widetilde{\chi}} ,
\endeq
where we take $\widetilde{\chi}=t-\left(\frac{1}{2}\widetilde{\nu}+\frac{1}{4}\right)\pi$.

Using (\ref{fip}) and Lemma \ref{lema-fiodt}, since function $\f(t)$ is continuous by the definition, there exists $k\in\mathbb{Z}$ such that
$$
\f(t)=(\widetilde{\chi}+k\pi)+O(t^{-1})\ \mathrm{ as\ } t\to\ty.
$$
For this $k$, we define $\bar{\phi}_0:=\tau_0+\left(k\pi-\left(\frac{1}{2}\widetilde{\nu}+\frac{1}{4}\right)\pi\right)$ and it holds
\begin{equation} \label{fiprekot}
\f(t)-\bar{\phi}_0=(t-\tau_0)+O(t^{-1})\ \textrm{as}\ t\to\ty.
\end{equation}
Define $\f_0:=\f(t_0)$ and notice that generally $\phi_0$ is not equal $\bar{\phi}_0$.

\smallskip

\emph{Step 3.} (Checking condition (\ref{biliplema_cond1}).) From (\ref{fiprekot}) it follows that $\f \simeq t\ \mathrm{as}\ t\to\infty$, and from (\ref{def_fazni_sustav}), using (\ref{bessel_J}--\ref{fun_S}) we get
\bgeq\label{rkvp}
r^2(t)=(x(t))^2+(-\dot x(t))^2 = \frac{2}{\pi} t^{-2+\mu} + O\left(t^{-3+\mu}\right)\ \mathrm{as}\ t\to\infty ,
\endeq
which implies
\bgeq\label{fkaotkaofi}
f(\f(t)) = r(t)\simeq t^{-\frac{2-\mu}{2}}\simeq \f^{-\frac{2-\mu}{2}}\ \mathrm{as}\ t\to\infty .
\endeq
Notice that from (\ref{def_fazni_sustav}) and (\ref{rkvp}) it follows that function $r(t)$ is of class $C^1$.

\smallskip

\emph{Step 4.} (Checking condition (\ref{biliplema_cond3}).)
By differentiating (\ref{fip}) we obtain
\bgeq\label{dfipodtp}
\frac{d\f}{dt}(t)=\cos ^2\f(t) \frac{\dot{x}(t)^2-x(t)\ddot{x}(t)}{x(t)^2}.
\endeq
Using (\ref{fip}) again, we have
\begin{equation} \label{cosnakvad_fiodt}
\cos^2 \f(t)=\frac{1}{1+\tan^2 \f(t)}= \frac{x(t)^2}{x(t)^2+\dot{x}(t)^2} .
\end{equation}
Substituting into (\ref{dfipodtp}) and using (\ref{bessel_J}--\ref{fun_S}) we get
\bgeq\label{dfi}
\lim_{t\to\infty}\frac{d\f}{dt}(t) = \lim_{t\to\infty}\frac{\dot{x}(t)^2-x(t)\ddot{x}(t)}{x(t)^2+\dot{x}(t)^2} = 1.
\endeq
On the other hand, differentiating (\ref{rkvp}) and from (\ref{def_fazni_sustav}), using (\ref{bessel_J}--\ref{fun_S}) again, as $-2+\mu<0$, we obtain that
\begin{eqnarray} \label{drpodtp}
\frac{dr}{dt}(t) & = & \frac{x(t) \dot{x}(t)+\dot{x}(t) \ddot{x}(t)}{\sqrt{x(t)^2+\dot{x}(t)^2}} = \frac{\frac{2}{\pi} t^{-3+\mu } (-2+\mu ) \sin(\widetilde{\chi})^2 + O\left(t^{-4+\mu}\right)}{r(t)} \nonumber\\
& \simeq & -t^{-\frac{2-\mu}{2}-1}\sin(\widetilde{\chi})^2+O\left(t^{-\frac{2-\mu}{2}-2}\right)\ \mathrm{as}\ t\to\infty .
\end{eqnarray}
Since $\frac{dr}{dt}(t)=f'(\f)\cdot\frac{d\f}{dt}(t)$ and since by (\ref{dfi}) we have  $\frac{d\f}{dt}(t)\simeq 1$ as $t\to\infty$, there exists $C_0, C_1\in\eR$, $C_1>C_0>0$, such that
\begin{equation} \label{dfleq}
|f'(\f)|\le C_0 t^{-\frac{2-\mu}{2}-1}\le C_1 \f^{-\frac{2-\mu}{2}-1}\ \textrm{as}\ \f\to\infty . \nonumber
\end{equation}

Also, notice that substituting (\ref{cosnakvad_fiodt}) in (\ref{dfipodtp}) it follows that function $\f(t)$ is of class $C^1$.

\smallskip

\emph{Step 5.} (Checking condition (\ref{biliplema_cond2b}).) From (\ref{dfi}) it follows that there exists $\tau_1\geq \tau_0$ such that $\frac{d\f}{dt}(t)>0$ for all $t\geq \tau_1$, so then function $\f(t)$ is strictly increasing for all $t$ sufficiently large. As $\f(t)$ is continuous, we conclude that for all $\f$ sufficiently large there exists inverse function $t=t(\f)$ of function $\f=\f(t)$ and that
$$
t(\f)= \f + \left(\tau_0-\overline{\phi}_0\right)+O(\f^{-1})\ \mathrm{as\ } \f\to\ty.
$$
Define $\phi_1:=\f(\tau_1)$ and notice that we can take $\tau_1$ sufficiently large such that $\phi_1\geq\phi_0$.

From (\ref{rkvp}), using more terms in the asymptotic expansion, we obtain
$$
r(t)=\sqrt{\frac{2}{\pi}}p_{\widetilde{\nu}}(t)\sqrt{\left(1+O\left(t^{-2}\right)\right) + \left(-\frac{p_{\widetilde{\nu}}'(t)}{p_{\widetilde{\nu}}(t)}+O\left(t^{-3}\right)\right)\sin 2\widetilde{\chi} + O\left(t^{-2}\right)\cos 2\widetilde{\chi}} ,
$$
where $p_{\widetilde{\nu}}(t)=t^{-\frac{1}{2}}P_{\widetilde{\nu}}(t)$. Using Lemma~\ref{tehnicka2} we conclude that for a fixed $\triangle\f>1$, we have
\begin{equation}\label{fiidelfi}
f(\f)-f(\f+\triangle\f)=r(t(\f))-r(t(\f+\triangle\f))\geq k_1\f^{-\frac{2-\mu}{2}-1},
\end{equation}
provided $\f$ is sufficiently large (to use Lemma~\ref{tehnicka2}, rescaling of $r(t)$ by factor $\sqrt{\frac{2}{\pi}}$ is needed). Moreover, by a careful examination of the proof of Lemma~\ref{tehnicka2}, we conclude that statement (\ref{fiidelfi}) uniformly holds for every $\triangle\f$ from a bounded interval, see Remark \ref{remark_prop_tehn2}, whose lower bound is greater than $1$, provided $\f$ is sufficiently large. But this means that before using Theorem \ref{biliplema}, we will have to ensure that $\theta$ from Theorem \ref{biliplema} is larger than $1$. It will be sufficient to take $\theta=(\pi/3+1)/2$, see \emph{Step 9}.

\smallskip

\emph{Step 6.} ($\Gamma'$ is a spiral near the origin.) Now we can prove that $\Gamma'$ is a spiral near the origin, that is, $f(\f)$ satisfies condition (\ref{def_spiral}) near the origin. First, from (\ref{fkaotkaofi}) it follows that $f(\f)\to 0$ as $\f\to\infty$. Second, from (\ref{fiidelfi}) it follows that $f(\f)$ is radially decreasing for all $\f$ sufficiently large, that is, there exists $\phi_2\geq\phi_1$ such that $f|_{[\phi_2,\infty)}$ is radially decreasing.

\smallskip

\emph{Step 7.} (The box dimension is invariant with respect to taking $\tau_0$ and $\phi_0$ sufficiently large.) First, we define $\tau_2$ to be such that $\f(\tau_2)=\phi_2$. Notice that $\tau_2$ is well-defined and $\tau_2\geq \tau_1$. From (\ref{rkvp}) it follows that $r(t)\geq 0$. Furthermore, $r(t)$ is strictly positive. Otherwise we would have $x(u_0)=0$ and $\dot{x}(u_0)=0$ for some $u_0\in\eR$ and from (\ref{gen_Bessel_equation}) would then follow that $x(u)=0$, for all $u\geq u_0$. This is a contradiction with $x(t)=\widetilde{J}_{\nu,\mu}(t)$ being a nontrivial solution of (\ref{gen_Bessel_equation}). As $r(t)$ is strictly positive, there exists constant $m_1>0$ such that for all $t\in[\tau_0,\tau_2]$ it holds
\begin{equation}\label{tvrdnja_vece}
r(t)>m_1 .
\end{equation}
Notice that $\phi_2\geq\phi_1\geq\phi_0$. From (\ref{fkaotkaofi}) it follows that $r(t)\to 0$ as $t\to\infty$ so there exists $\tau_3\geq\tau_2$ such that for all $t\in[\tau_3,\infty)$ it holds
\begin{equation}\label{tvrdnja_manje}
r(t)<m_1 .
\end{equation}
We define $\phi_3:=\f(\tau_3)$. Notice that we could increase $\tau_3$ and $\phi_3$ to accommodate all requirements in different parts of the proof on $t$ or $\f$ being sufficiently large. Now, from (\ref{tvrdnja_vece}) and (\ref{tvrdnja_manje}) we conclude that
\begin{equation}\label{komad_jedan}
\Gamma'|_{[\tau_0,\tau_2]}\bigcap\Gamma'|_{(\tau_3,\infty)}=\emptyset .
\end{equation}
As $f|_{[\phi_2,\infty)}$ is radially decreasing and function $\f(t)$ is strictly increasing for all $t\in[\tau_2,\infty)$, it follows that $\Gamma'|_{(\tau_2,\infty)}$ does not have self intersections, so
\begin{equation}\label{komad_dva}
\Gamma'|_{[\tau_2,\tau_3]}\bigcap\Gamma'|_{(\tau_3,\infty)}=\emptyset .
\end{equation}
Finally, from (\ref{komad_jedan}) and (\ref{komad_dva}) we have $\Gamma'|_{[\tau_0,\tau_3]}\bigcap\Gamma'|_{(\tau_3,\infty)}=\emptyset$. Now, we can apply Lemma \ref{novalemazuzu} on curve $\Gamma'$.

Using Lemma \ref{novalemazuzu} we see that without loss of generality we can assume that $\tau_0$ and $\phi_0$, in the assumptions of the theorem, are sufficiently large. Informally, we can always remove any rectifiable part from the beginning of $\Gamma'$, without changing the box dimension of $\Gamma'$.

\smallskip

\emph{Step 8.} (Checking waviness condition (\ref{def_SRC}).) By factoring (\ref{drpodtp}), again from (\ref{def_fazni_sustav}), using (\ref{bessel_J}--\ref{fun_Q}), we get
\begin{eqnarray}\label{drpodtfactp}
\frac{dr}{dt}(t) & = & \sqrt{\frac{2}{\pi}}\left( p_{\widetilde{\nu}}(t)+q_{\widetilde{\nu}}'(t)+\cot(\widetilde{\chi}) \left(q_{\widetilde{\nu}}(t)-p_{\widetilde{\nu}}'(t)\right) \right) \cdot \nonumber\\
& \cdot &  \left( 2 p_{\widetilde{\nu}}'(t)+q_{\widetilde{\nu}}''(t)+\cot(\widetilde{\chi}) \left(2 q_{\widetilde{\nu}}'(t)-p_{\widetilde{\nu}}''(t)\right) \right) \frac{\sin^2\widetilde{\chi}}{r(t)} ,
\end{eqnarray}
for every $\widetilde{\chi}\neq k\pi$, $k\in\Ze$ (division by factor $\sin\widetilde{\chi}$), where $p_{\widetilde{\nu}}(t)=t^{\frac{\mu}{2}-1}P_{\widetilde{\nu}}(t)$ and $q_{\widetilde{\nu}}(t)=t^{\frac{\mu}{2}-1}Q_{\widetilde{\nu}}(t)$.

Now, further using (\ref{fun_P}) and (\ref{fun_Q}), we get
\begin{eqnarray}
-\frac{p_{\widetilde{\nu}}(t)+q_{\widetilde{\nu}}'(t)}{q_{\widetilde{\nu}}(t)-p_{\widetilde{\nu}}'(t)} & \sim & -\frac{8}{8-6 \mu +\mu ^2+4 \nu ^2}\cdot t \quad\textrm{as}\ t\to\infty, \label{uvjet11_unique_lema} \\
-\frac{2p_{\widetilde{\nu}}'(t)+q_{\widetilde{\nu}}''(t)}{2q_{\widetilde{\nu}}'(t)-p_{\widetilde{\nu}}''(t)} & \sim & -\frac{8 (-2+\mu )}{(-4+\mu ) \left(4-4 \mu +\mu ^2+4 \nu ^2\right)}\cdot t \quad\textrm{as}\ t\to\infty. \label{uvjet12_unique_lema}
\end{eqnarray}

Analogously, we see that the derivatives
\begin{eqnarray}
\frac{d}{dt} \left(-\frac{p_{\widetilde{\nu}}(t)+q_{\widetilde{\nu}}'(t)} {q_{\widetilde{\nu}}(t)-p_{\widetilde{\nu}}'(t)} \right) & \sim & -\frac{8}{8-6 \mu +\mu ^2+4 \nu ^2} \quad\textrm{as}\ t\to\infty, \label{uvjet21_unique_lema} \\
\frac{d}{dt} \left(-\frac{2p_{\widetilde{\nu}}'(t)+q_{\widetilde{\nu}}''(t)} {2q_{\widetilde{\nu}}'(t)-p_{\widetilde{\nu}}''(t)} \right) & \sim & -\frac{8 (-2+\mu )}{(-4+\mu ) \left(4-4 \mu +\mu ^2+4 \nu ^2\right)} \nonumber\\ & & \quad\textrm{as}\ t\to\infty , \label{uvjet22_unique_lema}
\end{eqnarray}
are bounded for $t$ sufficiently large.

By Lemma \ref{unique} and using (\ref{uvjet11_unique_lema}--\ref{uvjet22_unique_lema}) (possibly with the shift by a constant, from variable $t$ to $\widetilde{\chi}$, which does not change asymptotic behavior), there exists $k_0\in\eN_0$ such that equations
\begin{equation} \label{tan_jedn}
       \cot\widetilde{\chi} = -\frac{p_{\widetilde{\nu}}(t)+q_{\widetilde{\nu}}'(t)}
       {q_{\widetilde{\nu}}(t)-p_{\widetilde{\nu}}'(t)} , \quad\quad
       \cot\widetilde{\chi} = -\frac{2p_{\widetilde{\nu}}'(t)+q_{\widetilde{\nu}}''(t)}
       {2q_{\widetilde{\nu}}'(t)-p_{\widetilde{\nu}}''(t)} ,
\end{equation}
have unique solutions $\hat{\chi}_{2k}$ and $\chi_{2k-1}$, respectively, in intervals $((k-1+k_0)\pi,(k+k_0)\pi)$, for each $k\in\eN_0$.

Notice that in (\ref{uvjet11_unique_lema}) and (\ref{uvjet12_unique_lema}) the right hand side coefficients multiplying variable $t$ are strictly less than zero, for every $\mu\in(0,2)$ and $\nu\in\eR$. So, from (\ref{uvjet11_unique_lema}) and (\ref{uvjet12_unique_lema}) follows that the right hand sides in (\ref{tan_jedn}) tend to minus infinity as $t\to\infty$. From the shape of function $\cot\widetilde{\chi}$ it follows that for every $\delta>0$ there exists $k_0(\delta)$ such that solutions  $\hat{\chi}_{2k}$ and $\chi_{2k-1}$ lie in intervals $((k+k_0(\delta))\pi-\delta,(k+k_0(\delta))\pi)$, for each $k\in\eN_0$. Here, we have possibly increased $k_0$ from Lemma \ref{unique} and re-indexed sequences $\hat{\chi}_{2k}$ and $\chi_{2k-1}$.

Also, notice that if $\nu\neq 0$ then the right hand side coefficient in (\ref{uvjet11_unique_lema}) is strictly less than the right hand side coefficient in (\ref{uvjet12_unique_lema}). So from (\ref{uvjet11_unique_lema}) and (\ref{uvjet12_unique_lema}) follows that if $\nu\neq 0$ then $\hat{\chi}_{2k}\neq \chi_{2k-1}$ for $k$ and analogously $t$ sufficiently large. Contrary, for $\nu=0$ the right hand side coefficients are equal. That case is treated in completely different manner and is discussed later at the end of the proof.

In terms of variables $\hat{t}_{2k}$ and $t_{2k-1}$ we see that for $\nu\neq 0$ and every $\delta>0$ there exists $k'_0(\delta)\in\eN_0$ (possibly again increased from the value of $k_0(\delta)$, to accommodate $k$ and $t$ being sufficiently large for $\hat{\chi}_{2k}\neq \chi_{2k-1}$ and sequences re-indexed) such that the equations (\ref{tan_jedn}) have unique and different solutions $\hat{t}_{2k}$ and $t_{2k-1}$, respectively, in intervals
\begin{equation} \label{interval_za_poz_t}
\left((k+k'_0(\delta))\pi+\left(\frac{1}{2}\widetilde{\nu}+\frac{1}{4}\right)\pi-\delta, (k+k'_0(\delta))\pi+\left(\frac{1}{2}\widetilde{\nu}+\frac{1}{4}\right)\pi\right) ,
\end{equation}
for each $k\in\eN_0$.

Without loss of generality we can take $t_{2k-1}<\hat{t}_{2k}$. Notice that $\hat{t}_{2k}-t_{2k-1}<\delta$. For this proof it will be sufficient to take $\delta=\pi/3$. It is easy to see from (\ref{drpodtfactp}) that $\frac{dr}{dt}(t)$ is positive between these solutions on intervals $(t_{2k-1},\hat{t}_{2k})$, for $t$ sufficiently large (equivalently taking $k'_0(\delta)$ to be sufficiently large).

As $\frac{d\f}{dt}(t)>0$ for all $t$ sufficiently large, from $\frac{dr}{dt}(t)=f'(\f)\cdot\frac{d\f}{dt}(t)$ it follows that $f'(\f)>0$ on set $\bigcup_{k=1}^{\ty}(\f_{2k-1},\hat{\f}_{2k})$ where $\f_{2k-1}:=\f(t_{2k-1})$ and $\hat{\f}_{2k}:=\f(\hat{t}_{2k})$.
This implies that function $f(\f)$ is increasing for some $\f$, so we can not apply Theorem \ref{novizuzu} directly.

Notice that if $t\in\bigcup_{k=0}^{\ty}(t_{2k-1},\hat{t}_{2k})$ then $r'(t)>0$ and if $t\in\bigcup_{k=0}^{\ty}(\hat{t}_{2k},t_{2k+1})$ then $r'(t)<0$.

We would like to prove that for every $k\in\eN_0$ there exists unique $t_{2k}\in(\hat{t}_{2k},t_{2k+1})$ such that $r(t_{2k})=r(t_{2k-1})$ and $t_{2k}-t_{2k-1}< \pi/3$ (where we will take $k'_0(\delta)$ from (\ref{interval_za_poz_t}) to be sufficiently large). As $r(\hat{t}_{2k})>r(t_{2k-1})$, and as function $r(t)$ is continuous and strictly decreasing on interval $(\hat{t}_{2k},t_{2k+1})$, it follows that, if such $t_{2k}$ exists then it is necessary unique, so we only need to prove the existence.

For every $k\in\eN_0$ we take $\bar{t}_{2k}:=t_{2k-1}+\pi/3$. Notice that $\bar{t}_{2k}\in(\hat{t}_{2k},t_{2k+1})$, because from (\ref{interval_za_poz_t}) follows that $t_{2k+1}-t_{2k-1}>\pi-\delta=2\pi/3$ and $\hat{t}_{2k}-t_{2k-1}<\delta=\pi/3$. Define $\bar{\f}_{2k}:=\f(\bar{t}_{2k})$ and take $\f_{2k-1}$ as defined before. Using (\ref{fiprekot}), we can take $t$ or equivalently $k'_0(\delta)$ sufficiently large, such that $(\pi/3+1)/2\leq\bar{\f}_{2k}-\f_{2k-1}\leq 2$ for every $k\in\eN_0$. (The exact value of the upper bound is not important. We just take some value larger than $\pi/3$. For lower bound, it is only important that it is larger than $1$ and lower than $\pi/3$, so we take the mean value between these two.)

Now, using Lemma~\ref{tehnicka2} and Remark \ref{remark_prop_tehn2}, analogously as in \emph{Step 5}, we compute
\begin{eqnarray}
r(t_{2k-1})-r(\bar{t}_{2k}) & = & r(t(\f_{2k-1}))-r(t(\bar{\f}_{2k})) \nonumber\\
& = & r(t(\f_{2k-1}))-r(t(\f_{2k-1}+(\bar{\f}_{2k}-\f_{2k-1}))) \nonumber\\
& \geq & C_2 \f_{2k-1}^{-\frac{2-\mu}{2}-1} > 0 , \nonumber
\end{eqnarray}
for some $C_2>0$, provided $\f$ or equivalently $k'_0(\delta)$ is sufficiently large. From this follows $r(\bar{t}_{2k})<r(t_{2k-1})$, and as function $r(t)$ is of class $C^1$, strictly decreasing on interval $(\hat{t}_{2k},\bar{t}_{2k})$ and $r(\hat{t}_{2k})>r(t_{2k-1})$, we see that there exist $t_{2k}\in(\hat{t}_{2k},\bar{t}_{2k})$ such that $r(t_{2k})=r(t_{2k-1})$ and obviously $t_{2k}-t_{2k-1}< \pi/3$.

Using $t_{2k+1}-t_{2k-1}>2\pi/3$, follows that $t_{2k+1}-t_{2k}>2\pi/3-\pi/3=\pi/3$.

We established that for every $k\in\eN_0$ holds $t_{2k+1}>t_{2k}>t_{2k-1}$. Notice that $r'(t_0)\leq 0$ and that sequence $(t_n)$, $n\in\eN_0$, is the same as the sequence from Definition \ref{def_wavy_function}, defined for function $r(t)$.

As $t_{2k+1}-t_{2k-1}>2\pi/3$ for every $k\in\eN_0$, we conclude that $t_n\to\infty$ as $n\to\infty$, which means that sequence $(t_n)$ satisfies condition (\ref{def_SRC})(i).

As $t_{2k+1}-t_{2k}>\pi/3$ for every $k\in\eN_0$, by taking $\varepsilon=\pi/3$, we see that sequence $(t_n)$ satisfies condition (\ref{def_SRC})(ii).

Using (\ref{drpodtp}) we conclude that there exist $C_3, C_4\in\eR$, $C_4>C_3>0$, such that
\begin{eqnarray}
\mathop{\mathrm{osc}}\limits_{t\in[t_{2k+1},t_{2k+2}]} r(t) & = & r(\hat{t}_{2k+2}) - r(t_{2k+1}) = \int\limits_{t_{2k+1}}^{\hat{t}_{2k+2}} r'(t)\,dt \nonumber\\
& \leq & \delta\cdot\sup\limits_{t\in[t_{2k+1},\hat{t}_{2k+2}]} r'(t) \leq C_3 t_{2k+1}^{-\frac{2-\mu}{2}-2} \leq C_4\hat{t}_{2k+2}^{-\frac{2-\mu}{2}-2} , \nonumber
\end{eqnarray}
for every $k\in\eN_0$, which means that sequence $(t_n)$ satisfies condition (\ref{def_SRC})(iii).

Finally, we conclude that, for $\nu\neq 0$, sequence $(t_n)$ satisfies waviness condition (\ref{def_SRC}), so $r(t)$ is a wavy function and $\Gamma'$ is a wavy spiral near the origin.

\smallskip

\emph{Step 9.} (Final conclusion.) From the previous steps we see directly that for $\nu\neq 0$, all of the assumptions of Theorem \ref{biliplema} are fulfilled. We take $\varepsilon'=(\pi/3+1)/2<\varepsilon$, $\theta=\min\{\varepsilon',\pi\}=(\pi/3+1)/2$ and $\a=\frac{2-\mu}{2}$. Using Theorem~\ref{biliplema} we prove that
\begin{equation} \label{th_final_concl}
\dim_B\Gamma'=\frac{4}{4-\mu} .
\end{equation}

\smallskip

\emph{Step 10.} (Degenerate case $\nu=0$.) Using known properties and transformations of Bessel functions, only in the case of $\nu=0$, expression (\ref{drpodtp}) for function $\dot{r}(t)$, from (\ref{def_fazni_sustav}), can be simplified in the form
$$
\frac{dr}{dt}(t) = \frac{t^{\mu-4}(\mu-2)}{4r(t)} \left(2tJ_{\frac{|\mu-1|}{2}-1}(t)+(\mu-1-|\mu-1|)J_{\frac{|\mu-1|}{2}}(t)\right)^2,
$$
which is less than or equal to $0$ for every $t\geq\tau_0$. This means that function $f(\f)$ is nonincreasing for $\f$ sufficiently large, so in this special case, we can use Theorem \ref{novizuzu} to prove claim (\ref{th_final_concl}).

\smallskip

The proof for $x(t)=\widetilde{Y}_{\nu,\mu}$ is analogous and will be omitted.

\qed

\appendix
\section{Auxiliary results}

\begin{prop}{\label{geometrijska}}
Let $XOY$ be an isosceles triangle such that $R=d(X,O)=d(Y,O)$. Let $\theta=\angle(XOY)\leq\pi$. Then
\begin{equation} \label{lema_geometrijska_claim}
\theta R \leq \frac{\pi}{2} d(X,Y) .
\end{equation}
\end{prop}

{\it Proof.}
Using elementary trigonometry we see that $d(X,Y)=2R\sin\frac{\theta}{2}$. Now, it is easy to see that
$$
\theta \leq \pi\sin\frac{\theta}{2},\quad\textrm{for all}\ \theta\in[0,\pi] ,
$$
which proves (\ref{lema_geometrijska_claim}).
\qed

\begin{lemma}{\label{lema-fiodt}}{\rm (Connection between $\f$ and $t$)}
Let $t_0>0$ and $\nu\in\eR$. Assume $\f:[t_0,\infty)\to\eR$ is a continuous function and let
\begin{equation}\label{tan-jedn}
\tan\f(t)=\frac{\sin(\chi)\left(1+O\left(t^{-2}\right)\right)+\cos(\chi)O\left(t^{-1}\right)}{\cos(\chi)\left(1+O\left(t^{-2}\right)\right)+\sin(\chi)O\left(t^{-1}\right)} \quad\textrm{as}\ t\to\infty ,
\end{equation}
where $\displaystyle \chi=t-\left(\frac{\nu}{2}+\frac{1}{4}\right)\pi$. Then there exists unique $k\in\Ze$ such that
\begin{equation} \label{lema-fiodt-rezultat}
\f(t)=(\chi+k\pi)+O\left(t^{-1}\right) \quad\textrm{as}\ t\to\infty .
\end{equation}

\end{lemma}

{\it Proof.}
Let $\chi$ is such that $|\sin\chi|\leq\frac{\sqrt{2}}{2}$. Then follows $|\cos\chi|\geq\frac{\sqrt{2}}{2}$ and $|\sin\chi|\leq|\cos\chi|\leq 1$. From this follows that the denominator of (\ref{tan-jedn}) equals $\cos(\chi)\left(1+O\left(t^{-1}\right)\right)$ so expression (\ref{tan-jedn}) becomes
$$
\frac{\sin(\chi)\left(1+O\left(t^{-2}\right)\right)+O\left(t^{-1}\right)}{\cos(\chi)\left(1+O\left(t^{-1}\right)\right)} \quad\textrm{as}\ t\to\infty .
$$
As lower and upper bounds on the denominator are strictly larger than zero and strictly lower then infinity, respectively, we write
$$
\frac{\sin(\chi)\left(1+O\left(t^{-2}\right)\right)}{\cos(\chi)\left(1+O\left(t^{-1}\right)\right)}+O\left(t^{-1}\right) = O\left(t^{-1}\right) + \tan(\chi)\left(1+O\left(t^{-1}\right)\right) \quad \textrm{as}\ t\to\infty .
$$
Using Lemma \ref{tehnicka0} and Lemma \ref{tehnicka1}, we compute
\begin{equation}\label{lema-fiodt-prvi-dio}
\tan\f(t)=\tan\left(\chi+O\left(t^{-1}\right)\right) ,
\end{equation}
for $|\sin\chi|\leq\frac{\sqrt{2}}{2}$ and $t$ sufficiently large.

Now let $\chi$ is such that $|\sin\chi|\geq\frac{\sqrt{2}}{2}$. Then follows $|\cos\chi|\leq\frac{\sqrt{2}}{2}$ and $1\geq |\sin\chi|\geq|\cos\chi|$. Analogously as before, we calculate the reciprocal expression of expression (\ref{tan-jedn}) to be equal to
\begin{eqnarray}
\cot\f(t) & = & \frac{\cos(\chi)\left(1+O\left(t^{-2}\right)\right)+\sin(\chi)O\left(t^{-1}\right)}{\sin(\chi)\left(1+O\left(t^{-2}\right)\right)+\cos(\chi)O\left(t^{-1}\right)} \nonumber \\
 & = &  O\left(t^{-1}\right) + \cot(\chi)\left(1+O\left(t^{-1}\right)\right) = \cot\left(\chi+O\left(t^{-1}\right)\right) \label{lema-fiodt-drugi-dio},
\end{eqnarray}
for $|\sin\chi|\geq\frac{\sqrt{2}}{2}$ and $t$ sufficiently large.

Finally, as function $\f(t)$ is a continuous function, from (\ref{lema-fiodt-prvi-dio}) and (\ref{lema-fiodt-drugi-dio}) it is easy to see that there exists unique $k\in\Ze$ such that it holds (\ref{lema-fiodt-rezultat}).

\qed

\begin{lemma}{\label{tehnicka0}}
For every $f:[t_0,\infty)\to\eR$, $t_0>0$, such that $\displaystyle\lim_{t\to\infty}f(t)=0$ there exist $y_i:[t_0,\infty)\to\eR$ such that $y_i(t)=O(f(t))$ as $t\to\infty$, $i=1,2$, and hold
\begin{enumerate}
\item[\rm{(i)}]  $(1+f(t))\tan t = \tan(t+y_1(t))$\quad for $t$ sufficiently large,
\item[\rm{(ii)}] $(1+f(t))\cot t = \cot(t+y_2(t))$\quad for $t$ sufficiently large.
\end{enumerate}
\end{lemma}

{\it Proof.}
\begin{enumerate}
\item[(i)]
For a fixed function $f(t)$ we pointwise show that there exists function $y_1(t)$ for $t$ sufficiently large. From the definition of the limit it follows that there exists $t_1\geq t_0$ such that $|f(t)|<\frac{1}{2}$, $\forall t\geq t_1$.

Now, let $t\geq t_1$. If $f(t)>0$ from Lemma \ref{tehnicka0A} follows that there exists $y_1(t)\in\eR$ such that $|y_1(t)|<\frac{\pi}{2}\cdot f(t)$ and it holds $(1+f(t))\tan t = \tan(t+y_1(t))$.
If $f(t)<0$, as $|f(t)|<\frac{1}{2}$, from Lemma \ref{tehnicka0B} follows that there exists value $y_1(t)\in\eR$ such that $|y_1(t)|<\pi\cdot (-f(t))$ and it holds $(1+f(t))\tan t = \tan(t+y_1(t))$.
Finally, if $f(t)=0$ the statement trivially holds for $y_1(t)=0$.

So, for constant $C=\pi$, hold $|y_1(t)|<C|f(t)|$ and $(1+f(t))\tan t = \tan(t+y_1(t))$ for $t$ sufficiently large.

\item[(ii)]
Define $f_1(t)=\frac{-f(t)}{1+f(t)}$. It is easy to see that $\displaystyle\lim_{t\to\infty}f_1(t)=0$ hence, using~(i), we compute
$$
(1+f(t))\cot t = \frac{1}{(1+f_1(t))\tan t} = \frac{1}{\tan(t+y_2(t))} = \cot(t+y_2(t)) ,
$$
for $t$ sufficiently large.
\end{enumerate}

\qed

\begin{lemma}{\label{tehnicka0A}}
Let $\a>0$. Then for every $t\in\eR$ there exists $y\in\eR$ such that $|y|<\frac{\pi}{2}\cdot \a$ and it holds $(1+\a)\tan t = \tan(t+y)$.
\end{lemma}
{\it Proof.}
Let $t\in\eR$. We distinguish four disjoint cases.

\textbf{Case 1}.
If $t=k\pi$ for some $k\in\Ze$. We take $y=0$, as it holds
$$
(1+\a)\tan k\pi = 0 = \tan(k\pi+y) .
$$

\textbf{Case 2}.
If $t=k\pi+\frac{\pi}{2}$ for some $k\in\Ze$. We also take $y=0$, because it holds
$$
(1+\a)\tan\left(k\pi+\frac{\pi}{2}\right)=\pm\infty=\tan\left(k\pi+\frac{\pi}{2}+y\right) .
$$

\textbf{Case 3}.
If $t\in(k\pi,k\pi+\frac{\pi}{2})$ for some $k\in\Ze$. We define $t_0:=t-k\pi$.

\begin{figure}
\begin{center}
\includegraphics{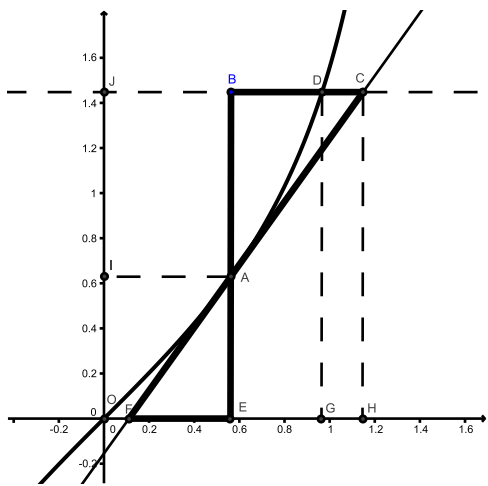}
\end{center}
\caption{Geometry in the proof of Lemma \ref{tehnicka0A}.}
\label{slika-tehnicka0A}
\end{figure}

Look at the graph of function $\tan t$ on interval $(k\pi,k\pi+\frac{\pi}{2})$, see figure \ref{slika-tehnicka0A}. Let point $E$ be on the $x$-axis corresponding to value $t$, point $I$ be on the $y$-axis corresponding to value $\tan t$ and point $A$ have coordinates $(t,\tan t)$. Let point $J$ on the $y$-axis correspond to $(1+\a)\tan t$.

As function $\tan t$ is a continuous and strictly increasing function on interval $(k\pi,k\pi+\frac{\pi}{2})$ and as statement $(1+\a)\tan t>\tan t$ holds, we conclude that there exists unique value $y$ such that it holds $(1+\a)\tan t=\tan(t+y)$. Let point $G$ on the $x$-axis correspond to value $t+y$. As function $\tan t$ is increasing function, notice that $y>0$ and $y=d(E,G)$. Also, let point $D$ have coordinates $(t+y,\tan(t+y))$.

Let point $C$ be the intersection of line $\overline{JD}$ parallel to the $x$-axis and the tangent on the graph of function $\tan t$ in point $A$. Because function $\tan t$ is convex function on interval $(k\pi,k\pi+\frac{\pi}{2})$, point $C$ is to the right of point $D$ and have coordinates $(t+y+z,\tan(t+y))$, where $z>0$ and $z=d(C,D)$. Let point $H$ be the projection of point $C$ to the $x$-axis, having coordinates $(t+y+z,0)$.

Next, let point $B$ be the intersection of the line $\overline{JD}$ parallel to the $x$-axis and the line $AE$ parallel to the $y$-axis. Notice that point $B$ have coordinates $(t,\tan(t+y))$.

Finally, let point $F$ be the intersection of the tangent on the graph of function $\tan t$ in point $A$ and the $x$-axis. Point $F$ is to the left of point $E$ and have the coordinates $(t-x,0)$, where $x>0$ and $x=d(E,F)$.

Now look at similar triangles $AEF$ and $ABC$. From similarity it follows
\begin{equation}\label{simequ}
\frac{d(A,E)}{d(E,F)}=\frac{d(A,B)}{d(B,C)} .
\end{equation}
From (\ref{simequ}) we compute
$$
\frac{x}{\tan t}=\frac{y+z}{(1+\a)\tan t - \tan t} > \frac{y}{\a\tan t} .
$$
Finally we get that $y<\a x$ and as function $\tan t$ is a convex function on interval $(k\pi,k\pi+\frac{\pi}{2})$, it follows that $x<t_0<\frac{\pi}{2}$ hence $y<\frac{\pi}{2}\cdot \a$. As $y>0$ it holds $|y|<\frac{\pi}{2}\cdot \a$.

\textbf{Case 4}.
If $t\in(k\pi+\frac{\pi}{2},(k+1)\pi)$ for some $k\in\Ze$. Now, $-t\in(-(k+1)\pi,-(k+1)\pi+\frac{\pi}{2})$. Using \textbf{Case 3}, we see that
$$
(1+\a)\tan t=-(1+\a)\tan(-t)=-\tan(-t+y)=\tan(t+(-y))
$$
and that $|-y|=|y|<\frac{\pi}{2}\cdot \a$.

\qed

\begin{lemma}{\label{tehnicka0B}}
Let $\a\in(0,\frac{1}{2})$. Then for every $t\in\eR$ there exists $y\in\eR$ such that $|y|<\pi\cdot \a$ and it holds $(1-\a)\tan t = \tan(t-y)$.
\end{lemma}
{\it Proof.} Let $t\in\eR$. We distinguish four disjoint cases.

\textbf{Case 1} and \textbf{Case 2}.
If $t=k\pi$ or $t=k\pi+\frac{\pi}{2}$, for some $k\in\Ze$, then we take $y=0$. The proof is analogous as in Lemma \ref{tehnicka0A}.

\textbf{Case 3}.
If $t\in(k\pi,k\pi+\frac{\pi}{2})$ for some $k\in\Ze$.

Look at the graph of function $\tan t$ on interval $(k\pi,k\pi+\frac{\pi}{2})$, see Figure~\ref{slika-tehnicka0A}.

Points $A$ to $J$ are defined analogously as in the proof of Lemma \ref{tehnicka0A}. The important difference here is that point $A$ has coordinates $(t-y,(1-\a)\tan t)$ and point $D$ has coordinates $(t,\tan t)$.

From similarity equation (\ref{simequ}), we compute
$$
\frac{x}{(1-\a)\tan t}=\frac{y+z}{\tan t - (1-\a)\tan t} > \frac{y}{\a\tan t} .
$$
We get that $y<\frac{\a}{1-\a}x<2\a x$ and as $x<\frac{\pi}{2}$ we get $y<\pi\cdot \a$. As $y>0$ it holds $|y|<\pi\cdot \a$.

\textbf{Case 4}.
If $t\in(k\pi+\frac{\pi}{2},(k+1)\pi)$ for some $k\in\Ze$. Analogous to the proof of Lemma \ref{tehnicka0A}.

\qed

\begin{lemma}{\label{tehnicka1}}
For every $f:\eR\to\eR$ there exist $y_i:\eR\to\eR$ such that $y_i(t)=O(f(t))$ as $t\to\infty$, $i=1,2$, and hold
\begin{enumerate}
\item[\rm{(i)}]  $f(t)+\tan t = \tan(t+y_1(t))\quad\textrm{for all}\ t\in\eR$,
\item[\rm{(ii)}] $f(t)+\cot t = \cot(t+y_2(t))\quad\textrm{for all}\ t\in\eR$.
\end{enumerate}
Moreover, function $y_1(t)$ can be chosen such that $|y_i(t)|\leq|f(t)|$ for all $t\in\eR$, $i=1,2$.
\end{lemma}

{\it Proof.}
\begin{enumerate}
\item[(i)]
For a fixed function $f(t)$ we will pointwise construct function $y_1(t)$, by defining the value $y_1(t)$ for every $t\in\eR$.

Let $t\in\eR$. We distinguish two disjoint cases.

\textbf{Case 1.}
If $t=k\pi+\frac{\pi}{2}$ for some $k\in\Ze$. We define $y_1(t):=0$, as for any value $f(t)\in\eR$, it holds
$$
f(t)+\tan t = \pm\infty =\tan(t) = \tan(t+y_1(t)) .
$$

\textbf{Case 2.}
If $t\in(k\pi-\frac{\pi}{2},k\pi+\frac{\pi}{2})$ for some $k\in\Ze$. Notice that function $\tan(t)$ is a continuous and strictly increasing function on interval $(k\pi-\frac{\pi}{2},k\pi+\frac{\pi}{2})$ and its image is set $\eR$. Therefore there exists unique value $t_1\in(k\pi-\frac{\pi}{2},k\pi+\frac{\pi}{2})$ such that $\tan(t_1)=\tan t +f(t)\in\eR$. We define $y_1(t):=t_1-t$, so we can write $\tan(t+y_1(t))=\tan t + f(t)$. By the mean value theorem we get
$$
f(t)=\tan(t+y_1(t))-\tan t = \tan'(\xi)\cdot y_1(t) ,
$$
for some value $\xi\in(k\pi-\frac{\pi}{2},k\pi+\frac{\pi}{2})$. As $\tan' (\xi)\geq 1$, for every $\xi\in(k\pi-\frac{\pi}{2},k\pi+\frac{\pi}{2})$ it follows that $|y_1(t)|\leq|f(t)|$.

\item[(ii)] Analogous as the proof of (i).

\end{enumerate}

\qed

\begin{lemma}\label{unique}
Let $F\in C^1(0,\ty)$ be such that $F(z)\sim az$ as $z\to\ty$ for some $a<0$. Assume that $\inf F'>-\ty$.
Then there exists a nonnegative integer $k_0$ such that for each $k\ge k_0$ the equation $\cot z=F(z)$ possesses the unique solution in
$J_k=(k\pi,(k+1)\pi)$.
\end{lemma}

The proof is given in \cite[Lemma 5]{cswavy}.

\bigskip

\bibliographystyle{abbrv}
\bibliography{bibliografija}

\end{document}